\documentclass[10pt]{article}
\usepackage{latexsym,amsmath,amssymb,graphics,amscd}
\usepackage{multirow}
\usepackage{booktabs}
\usepackage{tabularx}
\usepackage[hang,footnotesize]{caption}
\usepackage[pdftex]{graphicx}
\usepackage{graphicx}

\addtocounter{footnote}{1}
\makeatletter
\@addtoreset{table}{bsection}
\def\thetable{\thesection.\@arabic\c@table}
\def\fps@table{h, t}
\@addtoreset{equation}{section}

\makeatother

\newtheorem{theorem}{Theorem}[section]
\newtheorem{definition}[theorem]{Definition}
\newtheorem{lemma}[theorem]{Lemma}
\newtheorem{remark}[theorem]{Remark}
\newtheorem{proposition}[theorem]{Proposition}

\newtheorem{example}[theorem]{Example}

\newcommand{\bfi}{\bfseries\itshape}
\newsavebox{\savepar}

\makeatletter

\usepackage{paralist}

\begin{document}

\title{\textbf{Finite sample forecasting with estimated temporally aggregated linear processes}}
\author{Lyudmila Grigoryeva$^{1}$ and Juan-Pablo Ortega$^{2}$}
\date{}
\maketitle

\begin{abstract}
We propose a finite sample based predictor for estimated linear one dimensional time series models and compute the associated total forecasting error. The expression for the error that we present takes into account the estimation error. Unlike existing solutions in the literature, our formulas  require neither assumptions on the second order stationarity of the sample nor Monte Carlo simulations for their evaluation. This result is used to prove the pertinence of a new hybrid scheme that we put forward for the forecast of linear temporal aggregates. This novel strategy consists of carrying out the parameter estimation based on disaggregated data and the prediction based on the corresponding aggregated model and data. We show that in some instances this scheme has a better performance than the ``all-disaggregated'' approach presented as optimal in the literature. 
\end{abstract}

\bigskip


\makeatletter
\addtocounter{footnote}{1} \footnotetext{%
Corresponding author. D\'{e}partement de Math\'{e}%
matiques de Besan\c{c}on, Universit\'{e} de Franche-Comt\'{e}, UFR des
Sciences et Techniques. 16, route de Gray. F-25030 Besan\c{c}on cedex.
France. {\texttt{Lyudmyla.Grygoryeva@univ-fcomte.fr} }}
\makeatother
\makeatletter
\addtocounter{footnote}{1} \footnotetext{%
Centre National de la Recherche Scientifique, D\'{e}partement de Math\'{e}%
matiques de Besan\c{c}on, Universit\'{e} de Franche-Comt\'{e}, UFR des
Sciences et Techniques. 16, route de Gray. F-25030 Besan\c{c}on cedex.
France. {\texttt{Juan-Pablo.Ortega@univ-fcomte.fr} }}
\makeatother

\section{Introduction}

The success of parametric time series models as a tool of choice in many research fields is due in part to their good performance when it comes to empirical forecasting based on historical samples. Once a data generating process (DGP) has been selected and estimated for the forecasting problem at hand, there is a variety of well studied forecasting procedures and algorithms available in the literature. The most widespread loss function used in the construction of predictors is the mean square forecasting error (MSFE); see the monographs~\cite{BrocDavisYellowBook, Brockwell2002, MR1278033, luetkepohl:book} and references therein for detailed presentations of the available MSFE minimization-based techniques. This is the approach to prediction that we follow in this work; the reader is referred to~\cite{granger1969} or Section 4.2 in~\cite{granger:newbold} for forecasting techniques based on other optimality criteria.

The stochastic nature of the time series models that we consider implies that the forecasts produced with them, carry in their wake an error that cannot be minimized even if the parameters of the model are known with total precision; we refer to this as the {\bfi  characteristic error} of the model. Additionally,  all that is known in most applications is a historical sample of the variable that needs to be forecasted, out of which a model needs to be selected and estimated. There are well-known techniques to implement this, which are also stochastic in nature and that increase the total error committed when computing a forecast; we talk in that case of  {\bfi  model selection error} and {\bfi  estimation error}. All these errors that one incurs in at the time of carrying out a forecasting task are of different nature and much effort has been dedicated in the literature in order to quantify them in the case of linear multivariate VARMA processes. 

Most results obtained in this direction have to do with the combination of the estimation and the characteristic errors; this compound error is always studied assuming independence between the realizations of the model that are used for estimation and the ones used for prediction; we refer the reader to~\cite{Baillie79, Reinsel1980, Yamamoto1980, Yamamoto1981, Dufour1985, Basu1986, Lutkepohl1987, Samaranayake1988}. Explicit expressions for these errors in the VARMA context  are available in the monograph~\cite{luetkepohl:book}. Indeed, if we assume that the sample out of which we want to forecast is a realization of the unique stationary solution of a VAR model,  this error can be written down~\cite[page 97]{luetkepohl:book} using the time-independent autocovariance of the process; the situation in the VARMA context is more complicated and the expression provided~\cite[page 490]{luetkepohl:book} requires Monte Carlo simulations for its estimation. 

The knowledge regarding the error associated to model selection is much more rudimentary and research in this subject seems to be in a more primitive state. A good description of the state of the art  can be found in~\cite[page 318]{LutkepohlHandbook2006} as well as in~\cite[page 89]{Lutkepohl1986}, and references therein. We do not consider this source of forecasting error in our work and hence in the sequel we will use the denomination {\bfi  total error} to refer to the combination of the characteristic with the estimation errors.

In this paper we concentrate on one dimensional linear processes, a subclass of which is the ARMA family. The first contribution in this paper is the formulation of a MSFE based predictor that takes as ingredients a finite sample and the coefficients of a linear model estimated on it, as well as the computation of the corresponding total error. The main improvements that we provide with respect to preexisting work on this question are:
\begin{itemize}
\item We make no hypothesis on the second order stationarity of the sample at hand; in other words, we do not assume that the sample is a realization of the stationary solution of the recursions that define the model. Such a hypothesis is extremely difficult to test in small and finite sample contexts and it is hence of much interest to be able to avoid it.
\item The expression for the total forecasting error is completely explicit and does not require the use of Monte Carlo simulations.
\end{itemize}
The interplay between the characteristic error, the estimation errors, and the forecasting horizon is highly nonlinear and can produce surprising phenomena. For example, as it is well known, the characteristic error is an increasing function of the horizon, that is, the further into the future we forecast, the more error we are likely to commit. When we take into account the estimation error,  the total error may decrease with the forecast horizon! We study this finite sample related phenomenon with the total error formula introduced in Theorem~\ref{theor1} and illustrate it with an example in Section~\ref{Stock aggregation examples}.

The characterization of the total forecasting error that we described serves as a basis for the second main theme of this paper, namely, the interplay between multistep forecasting, the prediction of temporal aggregates, and the use of temporal aggregation estimation based techniques to lower the total forecasting error. In this part of the paper we work strictly in the ARMA context. The temporal aggregation of ARMA processes is a venerable topic that is by now well understood~\cite{Amemiya1972, Tiao1972, Brewer1973, TiaoWei1976, WeiTempAggr1979, Stram1986, Wei:book, Silvestrini2008} and has been extensively studied and exploited in the context of forecasting~\cite{abraham:aggregation, Lutkepohl1986, Lutkepohl1986a, Lutkepohl1987, Lutkepohl1989, Lutkepohl1989statPapers, Rossana1995, LutkepohlHandbook2006, Lutkepohl2009, Lutkepohl2010b} mainly by H. L\"utkepohl.

A recurrent question in this setup consists of determining the most efficient way to compute multistep forecasts or, more generally, predictions of linear temporal aggregates of a given time series. More specifically, given a sample and an underlying model, we can imagine at least two ways to construct a $h$ time steps ahead forecast,  or in general the one that is  a linear combination of the $h$ steps ahead values for the time series. First, we can simply compute the $h$ time steps ahead forecasts of the time series out of the original disaggregated sample and to determine the needed aggregate prediction out of them;  another possibility would be to temporally aggregate the sample and the time series model in such a way that the required forecast becomes a one time step ahead forecast for the new aggregated sample and model. If we do not take into consideration estimation errors and we only consider the characteristic error, there is a general result that states that the forecast based on high frequency disaggregate data has an associated error that is smaller or equal than the one associated to the  aggregate sample and model (we will recall it in Proposition~\ref{comparison of frequency predictors}). In the VARMA context, H. L\"utkepohl~\cite{Lutkepohl1986, Lutkepohl1987, Lutkepohl2009} has characterized the situations in which there is no loss of forecasting efficiency when working with temporally aggregated ingredients.

When estimation errors are taken into account, the inequality that we just described becomes strict~\cite{Lutkepohl1986, Lutkepohl1987}, that is, forecasts based on models estimated using the dissagregated high-frequency samples perform always better than those based on models estimated using aggregated data. This is so even in the situations described in~\cite{Lutkepohl1986, Lutkepohl1987, Lutkepohl2009} for which the characteristic errors associated to the use of the aggregated and the disaggregated models are identical; this is intuitively very reasonable due to the smaller sample size associated to the aggregated situation, which automatically causes an increase in the estimation error.

In Section~\ref{A hybrid forecasting scheme using aggregated time series models} we propose a forecasting scheme that is a hybrid between the two strategies that we just described. We first use the high frequency data for estimating a model. Then, we temporally aggregate the data and the model and finally forecasting is carried out based on these two aggregated ingredients. 
We will show that this scheme presents two major advantages:
\begin{itemize}
\item The model parameters are estimated using all the information available with the bigger sample size provided by the disaggregated data. Moreover, these parameters can be updated as new high frequency data becomes available.
\item In some situations, the total error committed using this hybrid forecasting scheme is {\it smaller} than the one associated to the forecast based on the disaggregated data and model and hence our strategy becomes optimal. Examples in this direction for both stock and flow temporal aggregates are presented in Section~\ref{Numerical results}. The increase in performance obtained with our method comes from minimizing the estimation error; given that the contribution of this error to the total one for univariate time series models is usually small for sizeable samples, the differences in forecasting performance that we will observe in practice are moderate. As we will show in a forthcoming work, this is likely to be different in the multivariate setup where in many cases, {\it the estimation error is  the main source of error}.
\end{itemize}

To our knowledge, this forecasting scheme has not been previously investigated in the literature and the improvement stated in the last point seems to be the first substantial application of temporal aggregation techniques in the enhancing of forecasting efficiency.

\section{Finite sample forecasting of linear processes}
In this section we introduce notations and definitions used throughout the paper and describe the framework in which we work. Additionally, since we are interested in finite sample based forecasting, we spell out in detail the predictors as well as the information sets on which our constructions are based.

\subsection{Linear processes}
Let $ \varepsilon = \left\{ \varepsilon _t \right\} _{ t=-\infty} ^{ \infty} $ be a set of independent and identically distributed random variables with mean zero and variance $ \sigma ^2 $. We will write in short 
\begin{equation*}
\varepsilon = \left\{ \varepsilon _t \right\}  \sim \text{IID} \left( 0, \sigma ^2 \right).
\end{equation*}
We say that $ X = \left\{ X _t \right\} _{ t=-\infty} ^{ \infty} $ is a {\bfi  linear causal} process whenever it can be represented as
\begin{equation} \label{LCP}
X _t = \sum^{\infty}_{i = 0} \psi _i \varepsilon _{ t-i}, \enspace \text{for all} \enspace t \in \mathbb{Z},
\end{equation}
where $ \left\{ \psi _i \right\} _{ i=0} ^{ \infty} $ is a set of real constants such that $ \sum^{\infty}_{i= 0} |{\psi _i }|<\infty $. Expression~\eqref{LCP} can also be rewritten as 
\begin{equation*}
X =\mathbf{ \Psi} \left( L \right) \varepsilon,
\end{equation*}
where $ L $ is the backward shift operator and $ \mathbf{\Psi} \left( z \right) $ is the power series $ \mathbf{\Psi} \left( z \right) = \sum^{\infty}_{i = 0} \psi _i z ^i  $.
The process $ X $ defined in \eqref{LCP} is called {\bfi  invertible} if there exist constants $ \left\{ \pi _j \right\} _{ j = 0} ^{ \infty} $ such that $ \sum^{\infty}_{j = 0} | \pi _j | < \infty $ and 
\begin{equation}
\label{LIP}
\varepsilon_t = \sum^{\infty}_{j = 0}  \pi_j X_{t - j}, \enspace \text{for all} \enspace t \in \mathbb{Z},
\end{equation}  
or equivalently,
$\varepsilon = \boldsymbol{\Pi} \left( L \right) X$,
where $\boldsymbol{\Pi} \left( z \right) $ is the power series $ \boldsymbol{\Pi} \left( z \right) = \sum^{\infty}_{j = 0} \pi _j z ^j  $. $\mathbf{\Psi} \left( L \right) $ and $\boldsymbol{\Pi} \left( L \right) $ can also be referred to as causal linear filter and invertible linear filter, respectively.
\subsection{Finite sample forecasting of causal and invertible ARMA processes}

Consider the causal and invertible ARMA(p, q)  specification determined by the equivalent relations
\begin{equation}  
 \label{1}
\mathbf{\Phi} \left( L \right) X_t = \mathbf{\Theta} \left( L \right) \varepsilon_t,  \enspace 
X_t = \sum^{\infty}_{i = 0} \psi_i  \varepsilon_{t - i}, \enspace 
\varepsilon_t = \sum^{\infty}_{j = 0}  \pi_j X_{t - j}.
\end{equation}
The innovations process $\varepsilon = \left\{ \varepsilon _t\right\} $ can be either independent and identically distributed IID$ (0, \sigma^{2}) $ or white noise WN$ (0, \sigma^{2}) $. In this subsection we  focus on how to forecast out of a finite sample $ \xi_{T} = \left\{ x_{1}, ..., x_{T} \right\} $ that satisfies the relations \eqref{1} and that has been generated out of a presample  $\left\{ x _{1-p}, ..., x_{0} \right\} $ and a preinnovations set $ \left\{ \varepsilon _{1 - q}, ..., \varepsilon_{0}\right\} $. A standard way to solve this problem~\cite{BrocDavisYellowBook, Brockwell2002} consists of assuming that $ {\xi_T} $ is a realization of the unique stationary process $ X $ that satisfies the ARMA relations \eqref{1} and to use its corresponding time independent autocovariance functions to formulate a linear system of equations whose solution provides the linear projection $ \widehat{X_{T + h}} $ of the random variable $ X_{T + h} $  onto $ \left\{ x_{1}, ..., x_{T} \right\} $ using the $L ^2 $ norm; this projection $ \widehat{X_{T + h}} $ minimizes the mean square error. We recall that writing the unique stationary solution of \eqref{1} usually requires knowledge about the infinite past history of the process. For example, for an AR(1) model of the form $ X_{t} - \phi X _{t - 1} = \varepsilon_{t} $, the unique stationary solution is given by $ X_t = \sum^{\infty}_{i = 0} \phi^i \varepsilon_{t - i} $.

Given that we are concentrating in the finite sample context, we prefer for this reason to avoid the stationarity hypothesis and the use of the corresponding autocovariance functions and to exploit in the forecast only the information that is strictly available, that is:
\begin{description}
\item  [(i)]  The model specification \eqref{1}: we assume that the model parameters are known with certainty and we neglect  estimation errors.
\item [(ii)] The sample $ \xi_T = \left\{ x_{1}, ..., x_{T} \right\} $.
\item [(iii)] The presample $\left\{ x _{1 - p}, ..., x_{0} \right\} $ and preinnovations $ \left\{ \varepsilon _{1 - q}, ..., \varepsilon _{0}\right\} $ that have been used in the sample generation.
\end{description}
We now define the preset $ I $ as
\begin{equation*} 
I: = \left\{
\begin{array}{l l}
\left\{ x_{1 - p}, ..., x_{0} \right\} \cup \left\{ \varepsilon _{1 - q}, ..., \varepsilon _{0}\right\} , &\text{when }   p, q \neq 0 \\
\left\{ x_{1 - p}, ..., x_{0} \right\} , &\text{when }  q = 0 \\
\left\{ \varepsilon _{1 - q}, ..., \varepsilon _0\right\} , &\text{when } p = 0.
\end{array} \right.
\end{equation*} 
Let $ r = \max \left\{ p, q \right\} $ and define the enlarged preset $ I ^{*} $ as 
$$I^{*}:= \left\{ x_{1 - r}, ..., x_{0} \right\} \cup \left\{ \varepsilon_{1 - r}, ..., \varepsilon_{0}\right\},$$
where:
\begin{itemize}
\item 
if $ p>q $: $ r = p $  and $ \varepsilon_{t} := \sum^{t+p-1}_{j = 0}  \pi_j X_{t - j}, \enspace 1 - p \le t < 1 - q $;
\item 
if $ q>p $: $ r= q $ and $ x_{t} = \sum^{t+q-1}_{i = 0} \psi_i  \varepsilon_{t - i}, \enspace 1 - q \le t < 1 - p $;
\item 
if $ q = p $:  $ I = I^{*} $.
\end{itemize}
The enlarged preset $ I^{*} $ is defined as a function of the elements in $ I $. Consequently, the sigma-algebras  $\sigma ( I ) $ and $ \sigma ( I^{*} ) $ generated by $ I $ and $ I^{*} $, respectively, coincide, that is, 
$\sigma (I) = \sigma ( I^\ast  )$. 

The following result is basically known (see for example~\cite{MR1278033, luetkepohl:book}) but we include it in order to be explicit and self-contained about the forecasting scheme that we are using in the rest of the paper and also to spell out the peculiarities of the finite sample setup in which we are working. We include a brief proof in the appendix.

\begin{proposition} \label{prop1}
In the conditions that we just described:
\begin{enumerate}
\item[\bf{(i)}]
The information sets (sigma algebras) $ \sigma \left( \underline{\epsilon_{T}} \right) := \sigma \left( I, \varepsilon_{1}, ..., \varepsilon_{T} \right)$ and $\sigma \left( \underline{{\xi_{T}}} \right) := \sigma \left( I, x_{1}, ..., x_{T} \right) $ generated by the preset and the past histories of the innovations $ \epsilon_{T}:=\{ \varepsilon _1, \ldots, \varepsilon _T\} $ and the sample $\xi_T:=\{x _1, \ldots, x_T\} $ coincide, that is,
\begin{equation} \label{eq:prop11}
\sigma \left( {\xi_{T}} \right) = \sigma \left( {\epsilon_{T}} \right).
\end{equation} 
\item[\bf{(ii)}]
If the innovations process is IID (respectively WN) then the optimal multistep forecast $ \widehat{ X_{T+h}} $  (respectively optimal linear forecast) based on $\sigma \left({\xi_{T}}\right) $ that minimizes the mean square forecasting error (MSFE) $ E \left[ \left( X _{T+h} - \widehat{ X _{T+h}} \right) ^2 \right]   $  is:
\begin{equation} \label{eq:prop12}
\widehat{ X_{T + h}} = \sum^{T + h  - 1 + r}_{i = h}  \psi_{i} \varepsilon_{T + h - i} = \sum^{T - 1 + r}_{i = 0} \psi_{i + h} \varepsilon_{T - i} = \sum^{T - 1 + r}_{i = 0} \sum^{T - i - 1 + r}_{j = 0} \psi_{i + h}\pi_{j} X_{T - i-j}.
\end{equation} 
\item[\bf{(iii)}]
The MSFE associated to this forecast is
\begin{equation} 
\label{eq:prop13}
{\rm MSFE} \left( \widehat{ X_{T + h} }\right) = \sigma ^2 \sum^{h - 1}_{i = 0} {\psi_{i}} ^{2}.
\end{equation} 
\item[\bf{(iv)}]
For ARMA models, the forecasts constructed in \eqref{eq:prop12} for different horizons with respect to the same information set $ \mathcal{F}_{T} := \sigma \left( {\xi_{T}} \right) = \sigma \left( {\epsilon_{T}} \right) $, satisfy the following recursive formula:
\begin{equation} \label{eq:prop14}
\widehat{X_{T + h}} = \left\{
\begin{array}{l l}
\phi_{1} \widehat{X_{T + h - 1}} + ... + \phi_{p} \widehat{X_{T + h - p}} + \theta_{h} \varepsilon_{T} + ... + \theta_{q} \varepsilon_{T + h - q}, & q \ge h,\\
\phi_{1} \widehat{X_{T + h - 1}} + ... + \phi_{p} \widehat{X_{T + h - p}}, & q < h.\\
\end{array} \right.
\end{equation} 
\end{enumerate} 
\end{proposition}
\begin{remark}
\label{remark forecast}
\normalfont
Testing the stationarity of small or in general  finite samples is a difficult task in practice. We emphasize that the prediction in Proposition~\ref{prop1} does not require any stationarity hypothesis. Moreover, we underline that the forecast~(\ref{eq:prop12}) does not coincide in general neither with the standard linear forecast for second order stationary series that uses the corresponding time independent autocovariance function (see for example~\cite{BrocDavisYellowBook}, page 63), nor with the usual finite sample approximation to the optimal forecast (see~\cite{MR1278033}, page 85). The main difference with the latter consists of the fact that the innovations associated to the presample are not assumed to be equal to zero but they are reconstructed out of it so that there is no loss of information. In the examples~\ref{MA1example} and~\ref{ARMA11} below we show how our forecast allows us to construct a predictor that:
\begin{description} 
\item [(i)]
is different from the one obtained assuming stationarity;
\item [(ii)]
has a better performance in terms of characteristic forecasting error.
\end{description} 
These statements do not generalize to arbitrary ARMA models; for example, for pure AR models, the predictor that we propose and those cited above coincide. \quad $\blacksquare$
\end{remark}

\begin{example}{Finite sample forecasting for the {\rm MA(1)} process.}
\label{MA1example}
\normalfont\\
We consider the MA(1) model
\begin{equation} 
\label{MA1}
X _t = \varepsilon _t + \theta \varepsilon _{ t-1}
\end{equation} 
and the trivial sample consisting of just one value $ x _1 $ at time $ t = 1 $; this sample is generated by the preset $ I = \left\{ \varepsilon _0 \right\} $ and the innovation $ \varepsilon _1 $. In this case, the enlarged preset $ I ^* = \left\{ x _0 , \varepsilon _0 \right\} $ with $ x _0 = \varepsilon _0 $. \\
Moreover, we have
\begin{itemize}
\item
$ \psi _0 = 1, \enspace \psi _1 = \theta$ and $ \psi _i = 0$, for any integer $ i > 1$,  
\item
$ \pi_0 = 1, \enspace \pi _j = \left( -1 \right) ^j \theta ^j $, for any integer $ j \geq 1$.
\end{itemize}
Consequently by \eqref{eq:prop12}, the forecast $ \widehat{ X _2 }$ based on the information set $\mathcal{F} _1=\sigma \left(\{I, x _1\} \right)$,  is given by 
\begin{equation*} 
\widehat{ X _2 } = \theta \varepsilon _1 = \theta \left( x _1 - \theta x _0 \right) = \theta x _1 - \theta ^2 \varepsilon _0,
\end{equation*} 
and has the associated error
\begin{equation*} 
\text{MSFE}( \widehat{ X_2 } ) = \sigma ^2.
\end{equation*}
On the other hand, the forecast that assumes that $ x _1 $ is a realization of the unique stationary solution of \eqref{MA1} and that uses the corresponding autocovariance function~\cite[page 63]{BrocDavisYellowBook}  is given by 
\begin{equation*}
\widehat{ X _2 ^S } := \dfrac{\theta }{1 + \theta ^2 } x _1,
\end{equation*}  
and has the associated error
\begin{equation*}
\text{MSFE}( \widehat{ X _2 ^S } ) = \sigma ^2 \left( 1+ \theta ^2 \right) - \dfrac{\sigma ^2 \theta ^2}{1+ \theta ^2 }. 
\end{equation*} 
We note that 
\begin{equation*}
\text{MSFE} ( \widehat{ X _2 ^S } ) = \sigma ^2 \left[ \left( 1+ \theta ^2 \right) - \dfrac{\theta ^2 }{1+ \theta ^2 } \right] > \sigma ^2 = \text{MSFE} ( \widehat{ X _2} ),
\end{equation*} 
which shows that the forecast that we propose has a better performance than the one based on the stationarity hypothesis. \quad $\blacksquare$
\end{example}
The better performance of the forecast that we propose in the preceding example can be in part due to the fact that we are using for $\widehat{ X _2 }$ additional information on the preinnovations that is not taken advantage of at the time of writing $\widehat{ X_2 ^S } $. In the following example we consider an ARMA(1,1) model and we see that the statements of Remark~\ref{remark forecast}  also hold, even though in this case, unlike in the MA(1) situation, the information sets on which the two forecasts considered are based are identical. 

\begin{example}{Finite sample forecasting for the {\rm ARMA(1,1)} process.}
\label{ARMA11}
\normalfont\\
Consider the model
\begin{equation*} 
X_t - \phi X _{ t-1} = \varepsilon _t + \theta \varepsilon _{ t-1} .
\end{equation*}
Then,
\begin{itemize}
\item
$ \pi _0 = 1 $, $ \pi _j = \left( -1 \right) ^j \left( \phi + \theta \right) \theta ^{ j-1} $,  for any integer $ j \geq 1$,
\item
$ \psi _0 = 1 $, $ \psi _i = \left( \phi + \theta \right) \phi ^{ i-1} $, for any integer $ i \geq 1$.
\end{itemize}
We consider the trivial sample $ x _1 $ generated by the preset $ I = \left\{ x _0 , \varepsilon _0 \right\} = I^* $. Using  Proposition~\ref{prop1}, we have that the one-step ahead forecast $ \widehat{ X _2 }$ based on the information set $\mathcal{F} _1=\sigma \left(\{I, x _1\} \right)$,  is given by
\begin{equation*} 
\widehat{X _2 } = \left( \phi + \theta \right) x _1 - \theta \left( \phi + \theta \right) x _0, \quad \mbox{with} \quad\enspace \text{MSFE} ( \widehat{ X _2 }) = \sigma ^2.
\end{equation*}
On the other hand, the forecast based on the stationarity hypothesis using the same information set, yields 
\begin{equation*}
\widehat{ X _2 ^S} := \dfrac{\left( \theta ^2 + \phi \theta + 1 \right) \left( \theta + \phi \right) \left( \phi \theta + 1\right) }{(\theta ^2 +\theta \phi +  1 )^2- \theta  ^2} x _1 - \dfrac{\left( \theta + \phi \right) \left(\theta \phi + 1 \right) \theta }{(\theta ^2 + \theta \phi +  1 )^2- \theta  ^2 }x _0,
\end{equation*}
and
\begin{equation*}
\text{MSFE}  ( \widehat{  X_2 ^S} )  = \dfrac{\left( \theta ^2 + \phi \theta + 1 \right) \left(\theta ^4  + \theta ^3 \phi+ \theta \phi + 1\right) }{(\theta ^2 +\theta \phi +  1 )^2- \theta  ^2 } \sigma ^2.
\end{equation*}

It is easy to check that the statement
$
\text{MSFE} ( \widehat{ X _2 ^S }) > \text{MSFE} ( \widehat{ X_2 } ) $ is equivalent to $ \theta ^4 \left(  \theta + \phi \right) ^2 > 0$,
which is always satisfied and shows that the forecast that we propose has a better performance than the one based on the stationarity hypothesis. \quad $\blacksquare$
\end{example}

\section{Forecasting with estimated linear processes}
\label{Forecasting with estimated linear processes}

In Proposition~\ref{prop1} we studied forecasting when the parameters of the model are known with total precision. In this section we explore a more general situation in which the parameters are estimated out of a sample. Our goal is to quantify the joint mean square forecasting error that comes both from the stochastic nature of the model (characteristic error) and the estimation error; we will refer to this aggregation of errors as the total error. This problem has been extensively studied in the references cited in the introduction always using the following two main constituents:
\begin{itemize}
\item Estimation and the forecasting are carried out using independent realizations of the time series model.
\item The model parameter estimator is assumed to be asymptotically normal (for example, the maximum likelihood estimator); this hypothesis is combined with the use of the so called {\bfi  Delta Method}~\cite{Serfling1980} in order to come up with precise expressions for the total error.
\end{itemize}

The most detailed formulas for the total error in the VARMA context can be found in~\cite{luetkepohl:book} where the Delta Method is applied to the forecast considered as a smooth function of the model parameters. If we assume that the sample out of which we want to forecast is a realization of the unique stationary solution of a VAR model, an explicit expression for this error can be written down by following this approach~\cite[page 97]{luetkepohl:book} that involves the time-independent autocovariance of the process. In the VARMA setup, the situation is more complicated~\cite[page 490]{luetkepohl:book} and the resulting formula requires the use of a Monte Carlo estimation. 

In subsection~\ref{The total error of finite sample based forecasting} we start by obtaining a formula for the total error using a different approach at the time of invoking the Delta Method; our strategy uses this method at a more primitive level by considering the parameters of the linear representation of the process seen as a function of the ARMA coefficients. We show that discarding higher order terms on $1/\sqrt{T} $, where $T $ is the sample size used for estimation, the resulting formula for the total error can be approximated by a completely explicit expression that  involves only the model parameters and the covariance matrix associated to the asymptotically normal estimator of the ARMA coefficients.

In subsection~\ref{On Lutkepohls formula for the total forecasting error} we rederive the total error formula by H. L\"utkepohl~\cite[page 490]{luetkepohl:book} and show that it can be rewritten as explicitly as ours without using any stationarity hypothesis or Monte Carlo simulations. Moreover, we show that this formula coincides with the approximated one obtained in subsection~\ref{The total error of finite sample based forecasting} by discarding higher order terms on $1/\sqrt{T} $.

\subsection{The total error of finite sample based forecasting}
\label{The total error of finite sample based forecasting}
Consider the causal and invertible ARMA(p, q) process $ \left\{ X_t \right\} $ determined by the equivalent relations
\begin{equation}
\label{31}  
\mathbf{\Phi} \left( L \right) X_t = \mathbf{\Theta} \left( L \right) \varepsilon_t,  \enspace 
X_t = \sum^{\infty}_{i = 0} \psi_i  \varepsilon_{t - i}, \enspace 
\varepsilon_t = \sum^{\infty}_{j = 0}  \pi_j X_{t - j}, \enspace
\left\{ \varepsilon_t \right\} \sim \text{IID} ({0}, \sigma^{2}),
\end{equation} 
and denote $ \mathbf{\Psi} :=\{ \psi_0, \psi_1, \ldots\}$, $ \mathbf{\Pi} :=\{ \pi_0, \pi_1, \ldots\}$.
In Proposition~\ref{prop1} we studied forecasting for the process \eqref{31} when the parameters $ \mathbf{\Psi} $ or $\mathbf{\Pi}$ of the model are known with total precision; in this section we suppose that these parameters are estimated by using a sample independent from the one that will be used for forecasting. A more preferable assumption would have been that the parameters  $\mathbf{\Psi}$ are estimated based on the same sample that we intend to use for prediction, but exploiting only data up to the forecasting origin; Samaranayake~\cite{Samaranayake1988} and Basu {\it et al}~\cite{Basu1986}  have shown that many results obtained in the presence of  the independence hypothesis remain valid under this more reasonable assumption.

Under the independence hypothesis, the model coefficients  $ \mathbf{\Psi} $ or $ \mathbf{\Pi} $ become random variables independent from the process $ X$ and the innovations $ \varepsilon $. Moreover, we assume that these random variables are asymptotically normal, as for example in the case of maximum likelihood estimation of the ARMA coefficients. 

For the sake of completeness, we start by recalling the Delta Method, that will be used profusely in the following pages. A proof and related asymptotic statements can be found in~\cite{Serfling1980}.

\begin{lemma}[Delta Method]
\label{delta method}
Let $\widehat{\boldsymbol{\beta}} $  be an asymptotically normal estimator for the vector parameter $\boldsymbol{\beta} \in \mathbb{R}^n $, that is, there exists a covariance matrix $\Sigma $ such that
\begin{equation*}
\sqrt{T} \left(\widehat{\boldsymbol{\beta}}- \boldsymbol{\beta}\right)\xrightarrow[T \rightarrow \infty]{{\rm dist}}N(0, \Sigma),
\end{equation*}
where $T$ is the sample size. Let $f: \mathbb{R} ^n\rightarrow \mathbb{R} ^m $ be a vector valued continuously differentiable function and let $J _f  $ be its Jacobian matrix, that is, $\left((J _f)(\boldsymbol{\beta})\right)_{ij}:=(\partial f _i/\partial \beta _j) (\boldsymbol{\beta}) $. If $J _f(\boldsymbol{\beta})\neq 0 $, then
\begin{equation*}
\sqrt{T} \left(f(\widehat{\boldsymbol{\beta}})- f(\boldsymbol{\beta})\right)\xrightarrow[T \rightarrow \infty]{{\rm dist}}N(0, J _f\Sigma J _f ').
\end{equation*} 
\end{lemma}

The next ingredient needed in the formulation of the main result of this section is the covariance matrix $\Sigma_{\boldsymbol{\Xi}_P}$ associated to the asymptotic normal character of the estimator $ \widehat{ \boldsymbol{\Xi} }_P $ for the parameters $ \boldsymbol{\Xi} _P:= \left( \psi _1, \ldots, \psi_P, \pi _1, \ldots, \pi _P\right)  $, for some integer $P$. This is spelled out in the following lemma whose proof is a straightforward combination of the Delta Method with the results in Section 8.8 of~\cite{BrocDavisYellowBook}.

\begin{lemma}
\label{xis and so on}
Let $\{X_t\} $ be a causal and invertible ARMA(p,q) process like in~(\ref{31}). Let $ \boldsymbol{\Phi} :=( \phi _1, \ldots, \phi_p)' $,   $\boldsymbol{\Theta}:=( \theta_1, \ldots, \theta_q)' $, and $\boldsymbol{\beta}:=\left(\boldsymbol{\Phi}', \boldsymbol{\Theta}'\right)' $  be the ARMA parameter vectors and let $ \boldsymbol{\Xi} _P:= \left( \psi _1, \ldots, \psi_P, \pi _1, \ldots, \pi _P\right)  $ be a collection of length $2P $ of the parameters that provide the linear causal and invertible representations of that model. Then:
\begin{description}
\item [(i)] The maximum likelihood estimator $\widehat{\boldsymbol{\beta}} $ of $\boldsymbol{\beta} $ is asymptotically normal
\begin{equation*}
\sqrt{T} \left(\widehat{\boldsymbol{\beta}}- \boldsymbol{\beta}\right)\xrightarrow[]{{\rm dist}}N(0, \Sigma_{\boldsymbol{\beta}}), \quad \mbox{with} \quad \Sigma_{\boldsymbol{\beta}}= \sigma^2 
\left(
\begin{array}{cc}
E\left[\mathbf{U}_t \mathbf{U}_t'\right]&E\left[\mathbf{U}_t \mathbf{V}_t'\right]\\
E\left[\mathbf{V}_t \mathbf{U}_t'\right]&E\left[\mathbf{V}_t \mathbf{V}_t'\right]
\end{array}
\right)^{-1},
\end{equation*}
where $\mathbf{U}_t :=(U _t, \ldots, U_{t+1-p})' $, ${\bf V} _t :=(V _t, \ldots, V_{t+1-q})' $, and $\{U _t\} $ and $\{V _t\} $ are the autoregressive processes determined by 
\begin{equation*}
\mathbf{\Phi}(L)U _t= \varepsilon _t \quad \mbox{and } \quad \mathbf{\Theta}(L)V _t= \varepsilon _t.
\end{equation*}
\item [(ii)] Consider the elements in $\mathbf{\Xi}_P $ as a function of $\boldsymbol{\beta} $, that is, $\mathbf{\Xi}_P(\boldsymbol{\beta}) := \left( \psi _1(\boldsymbol{\beta}) , \dots, \psi_P(\boldsymbol{\beta}) , \pi _1 (\boldsymbol{\beta}), \dots, \pi _P(\boldsymbol{\beta}) \right) $. Then, by the Delta Method we have that
\begin{equation} 
\label{eq:theor33}
\sqrt{T} \left(\widehat{\mathbf{\Xi}}_P - \mathbf{\Xi}_P \right) \xrightarrow{\enspace d \enspace} N (0, \Sigma _{\mathbf{\Xi}_P} ),\quad \mbox{where} \quad \Sigma _{\mathbf{\Xi}_P}:=J _{\mathbf{\Xi}_P} \Sigma _{\mathbf{\boldsymbol{\beta}}} J_{\mathbf{\Xi}_P}'
\end{equation} 
and $ (J_{\mathbf{\Xi}_P})_{ij} = \dfrac{\partial{(\boldsymbol{\Xi}_P)_i}}{\partial{\beta _j}}$, $i=1, \ldots, 2P $,  $j=1,\dots,p+q$. Details on how to algorithmically compute the Jacobian $J_{\mathbf{\Xi}_P} $ are provided in Appendix~\ref{computation of the jacobian}.
\end{description}
\end{lemma}

The next theorem is the main result in this section. Its proof can be found in the appendix.

\begin{theorem}
\label{theor1}
Let $ \xi_T = \left\{ x_1, \dots , x_T\right\} $ be a sample obtained as a realization of the causal and invertible ARMA(p,q) model in~\eqref{31} using a preset $I$. In  order to forecast out of this sample, we first estimate the parameters of the model $ \widehat{\mathbf{\Psi}} = \left\{ \widehat{\psi} _0 ,  \widehat{\psi}_1, \dots \right\}  $, $ \widehat{\mathbf{\Pi}} = \left\{ \widehat{\pi} _0 ,  \widehat{\pi}_1, \dots \right\}  $ based on another sample $ \xi'_T$ that we assume to be independent of $\xi_T$, using the maximum likelihood estimator $\widehat{ \boldsymbol{\beta}}:= \left( \widehat{\boldsymbol{\Phi}}', \widehat{\boldsymbol{\Theta}}'\right)' $ of the ARMA parameters. If we use the forecasting scheme introduced in Proposition~\ref{prop1}, then:
\begin{enumerate}
\item[\bf{(i)}]
The optimal multistep forecast $\widehat{X_{T+h}} $ for $X_{T+h}$ based on the information set $\mathcal{F} _T $ generated by the sample $\xi_T$ and using the coefficients estimated on the independent sample $ \xi'_T$ is
\begin{equation} \label{eq:theor31}
\widehat{ X_{T + h}} = \sum^{T + h  - 1 + r}_{i = h} \widehat{ \psi}_{i} {\tilde{\varepsilon}}_{T + h - i},
\end{equation} 
where $ r= \max \{p,q \}$ and $ \tilde{\varepsilon }_t := \sum^{t+r-1}_{j = 0} \widehat{\pi}_j x_{t-j}$.
\item[\bf{(ii)}]
The mean square forecasting error (MSFE) associated to this forecast is
\begin{align} 
\label{eq:theor32}
 {\rm MSFE} \left( \widehat{ X_{T + h} }\right) &= \sigma ^2  \sum^{h - 1}_{i = 0} {\psi_{i}} ^{2} 
 + \sigma ^2 \Biggl[  \sum^{P}_{i = h} {\psi_{i}} ^{2} - 2 \sum^{P}_{i = h} \sum^{P-i}_{j = 0} \sum^{P-i-j}_{k = 0}  \psi_{i+j+k} \psi_k E \left[ \widehat{ \psi} _i  \widehat{\pi} _j  \right]  \nonumber \\
 & + \sum^{P}_{i = h} \sum^{P-i}_{j = 0} \sum^{P-i-j}_{k = 0} \sum^{P}_{i' = h} \sum^{P-i'}_{j' = 0} \sum^{P-i'-j'}_{k' = 0} \psi _k \psi _{k'} E \left[ \widehat{ \psi }_i  \widehat{\pi} _j \widehat{\psi}_{i'}
\widehat{\pi} _{j'} \right]  \delta _{i+j+k,i'+j'+k'}  
 \Biggr],
\end{align}
where $P = T+h-1+r$.

The first summand will be referred to as the \textbf{characteristic forecasting error} and the second one as the \textbf{estimation based forecasting error}. Notice that the characteristic error coincides with \eqref{eq:prop13} and amounts to the forecasting error committed when the model parameters are known with the total precision.
\item[\bf{(iii)}]
Let $\mathbf{\Xi}_P := \left( \psi _1 , \dots, \psi_P , \pi _1 , \dots, \pi _P \right) $ with $P = T+h-1+r $.   Using the notation introduced in Lemma~\ref{xis and so on} and discarding higher order terms in $1/\sqrt{T}$, the MSFE in~\eqref{eq:theor32} can be approximated by  
\begin{align} 
\label{eq:theor34}
 {\rm MSFE} \left( \widehat{ X_{T + h} }\right) &= \sigma ^2  \sum^{h - 1}_{i = 0} {\psi_{i}} ^{2} 
+ \sigma ^2 \dfrac{1}{T}\Biggl[  \sum^{P}_{i = h}(\Sigma _{\mathbf{\Xi}_P})_{i,i}  + 2 \sum^{P}_{i = h} \sum^{P-i}_{j = 0} \sum^{P-i-j}_{k = 0}\psi_{i} \psi_k  (\Sigma _{\mathbf{\Xi}_P})_{i+j+k,j+P} \nonumber \\
& + \sum^{P}_{i = h} \sum^{P-i}_{j = 0} \sum^{P-i-j}_{k = 0} \sum^{P}_{i' = h} \sum^{P-i'}_{j' = 0} \sum^{P-i'-j'}_{k' = 0} \psi _i \psi _k \psi _{i'} \psi _{k'}  (\Sigma _{\mathbf{\Xi}_P})_{j+P,j'+P} \delta _{i+j+k,i'+j'+k'} \Biggl],
\end{align} 
where $\Sigma_{\boldsymbol{\Xi}_P} $ is the covariance matrix in~(\ref{eq:theor33}).
\end{enumerate} 
\end{theorem}

\subsection{On L\"utkepohl's formula for the total forecasting error}
\label{On Lutkepohls formula for the total forecasting error}

As we already pointed out, H. L\"utkepohl~\cite[pages 97 and 490]{luetkepohl:book} has proposed formulas for VARMA models similar to the ones presented in Theorem~\ref{theor1} based on a different application of the Delta Method. In this section, we  rederive L\"utkepohl's result in the ARMA context and show that it is identical to the approximated formula~(\ref{eq:theor34}) presented in part {\bf (iii)} of Theorem~\ref{theor1}. In passing, this conveys that L\"utkepohl's result can be explicitly formulated and computed using neither stationarity hypotheses nor Monte Carlo simulations.

The key idea behind L\"utkepohl's formula is applying the Delta Method by thinking of the forecast $\widehat{ X_{T + h}} $ in question as a differentiable function $\widehat{ X_{T + h}}(\boldsymbol{\beta}) $ of the model parameters $\boldsymbol{\beta}:= (\mathbf{\Phi}',\mathbf{\Theta}')'$. In order to develop further this idea, consider first the information sets $\mathcal{F}_T:= \sigma(\underline{\xi_T})$ and $\mathcal{F}_T':= \sigma(\underline{\xi_T'})$ generated by two independent samples $\xi_T $ and $\xi _T' $ of the same size. The sample $\xi_T $ is used for forecasting and hence $\mathcal{F} _T $ determines the forecast $\widehat{ X_{T + h}}(\boldsymbol{\beta})$ once the model parameters $\boldsymbol{\beta} $ have been specified. 
The sample $\xi_T ' $ is in turn used for model estimation and hence $\mathcal{F}_T' $ determines $\widehat{\boldsymbol{\beta}} $. Consequently, the random variable $\widehat{ X_{T + h}}(\widehat{\boldsymbol{\beta}})$ is fully determined by $\mathcal{F} _T $ and $\mathcal{F}_T' $. In this setup,  
a straightforward application of the statement in Lemma~\ref{delta method} shows that
\begin{equation}
\label{delta use luetkepohl}
\sqrt{T} \left[ \widehat{X_{T+h}}\left( \widehat{\boldsymbol{\beta} }\right)  - \widehat{X_{T+h}} \left( \boldsymbol{\beta} \right)\mid \mathcal{F} _T \right]  \xrightarrow[T \rightarrow \infty]{\enspace dist \enspace} N \left( 0,  \dfrac{\partial{\widehat{X_{T+h}}}}{\partial{\boldsymbol{\beta} }} \Sigma _{\boldsymbol{\beta}} \left( \dfrac{\partial{\widehat{X_{T+h}}}}{\partial{\boldsymbol{\beta} }}  \right)' \right),
\end{equation}
which, as presented in the next result is enough to compute the total forecasting error.

\begin{theorem}
\label{luetkepohl formula theorem}
In the same setup as in Theorem~\ref{theor1}, the total error associated to the forecast in~(\ref{eq:theor31}) can be approximated by
\begin{equation}
\label{luetkepohl formula 1}
{\rm MSFE} \left( \widehat{ X_{T + h} }\right) = \sigma ^2  \sum^{h - 1}_{i = 0} {\psi_{i}} ^{2}  + \frac{1}{T}E \left[ \dfrac{\partial{\widehat{X_{T+h}}}}{\partial{\boldsymbol{\beta} }} \Sigma _{\boldsymbol{\beta}} \left( \dfrac{\partial{\widehat{X_{T+h}}}}{\partial{\boldsymbol{\beta} }}  \right)' \right].
\end{equation}
We refer to this expression as {\bfi  L\"utkepohl's formula} for the total forecasting error. Moreover:
\begin{description}
\item [(i)] Lutkepohl's formula coincides with the approximate expression for the total error stated in~(\ref{eq:theor34}). In particular, the second summand in L\"utkepohl's formula, which describes the contribution to the total error given by the estimation error, can be expressed as:
\begin{align}
\label{first expresion for luetkepohl}
\Omega(h) :=&E \left[ \dfrac{\partial{\widehat{X_{T+h}}}}{\partial{\boldsymbol{\beta} }} \Sigma _{\boldsymbol{\beta}} \left( \dfrac{\partial{\widehat{X_{T+h}}}}{\partial{\boldsymbol{\beta} }}  \right)' \right]=\sigma ^2 \Biggl[  \sum^{P}_{i = h}(\Sigma _{\mathbf{\Xi}_P})_{i,i}  
	+ 2 \sum^{P}_{i = h} \sum^{P-i}_{j = 0} \sum^{P-i-j}_{k = 0}\psi_{i} \psi_k  (\Sigma _{\mathbf{\Xi}_P})_{i+j+k,j+P} \nonumber \\ 
	& + \sum^{P}_{i = h} \sum^{P-i}_{j = 0} \sum^{P-i-j}_{k = 0} \sum^{P}_{i' = h} \sum^{P-i'}_{j' = 0} \sum^{P-i'-j'}_{k' = 0} \psi _i \psi _k \psi _{i'} \psi _{k'}  (\Sigma _{\mathbf{\Xi}_P})_{j+P,j'+P} \delta _{i+j+k,i'+j'+k'} \Biggl].
\end{align}
\item [(ii)] If we assume that the samples used for forecasting are second order stationary realizations of the model~(\ref{31}) and $\gamma: \Bbb Z \rightarrow \mathbb{R} $ is the corresponding time invariant autocovariance function, then the estimation error can be expressed as:
\begin{multline}
\label{second expresion for luetkepohl}
\frac{1}{T}\Omega(h)\\
=\frac{1}{T}\sum_{i=h}^P\sum_{j=0}^{P-i}\sum_{i'=h}^P\sum_{j'=0}^{P-i'}\left[\pi_j\pi_j'
(\Sigma _{\mathbf{\Xi}_P})_{i,i'}+2 \pi_j \psi_i' (\Sigma _{\mathbf{\Xi}_P})_{i,j'+P}+\psi_i\psi_i' (\Sigma _{\mathbf{\Xi}_P})_{j+P,j'+P}
\right]\gamma(i+j-i'-j').
\end{multline}
\end{description}
\end{theorem}

\section{Finite sample forecasting of temporally aggregated linear processes}
\label{Finite sample forecasting of temporally aggregated linear processes}

The goal of this section is proposing a forecasting scheme for temporal aggregates based on using  high frequency data for estimation purposes and the corresponding temporally aggregated model and data for the forecasting task. We show, using the formulas introduced in the previous section, that in some occasions this strategy can yield forecasts of superior quality  than those based exclusively on high frequency data that are presented in the literature as the best performing option~\cite{Lutkepohl1986, Lutkepohl1987, Lutkepohl2009}.

We start  by recalling general statements about temporal aggregation that we need in the sequel. We then proceed by using various extensions of the results in Section~\ref{Forecasting with estimated linear processes} regarding the computation of total forecasting errors with estimated series in order to compare the performances of the schemes that we just indicated.

\subsection{Temporal aggregation of time series}

The linear temporal aggregation of  time series requires the use of the elements provided in the following definition.

\begin{definition} \label{defTA}
Given $ K \in \mathbb{N} $, $ X $ a time series, and $ \mathbf{w} = \left( w _1 , ..., w _K \right) ' \in \mathbb{R}^K  $, define the $ K ${\bfi  -period projection} $ p _K $  of $X $ as 
\begin{equation*}
p _K \left( X \right) := \left\{ \mathbf{X _j ^{ \left( K \right)} } \right\} _{ j \in \mathbb{Z} } \in \prod_{j \in \mathbb{Z}} \mathbb{R}^K,
\end{equation*}
where $ \mathbf{X _j ^{(K)}} := \left( X_{(j-i) K +1}, ..., X _{j K}  \right) \in \mathbb{R}^{ K}, $ and the corresponding {\bfi  temporally aggregated time series} $ Y $ as 
\begin{equation} \label{TAseries}
Y: = \mathbb{I} _\mathbf{w} \circ p _K \left( X \right), 
\end{equation}
where \\
$$ \mathbb{I} _\mathbf{w}: \enspace \prod_{ j \in \mathbb{Z}} \mathbb{R} ^K \longrightarrow \prod _{ j \in \mathbb{Z}} \mathbb{R} , \enspace
\prod_{ j \in \mathbb{Z}} \mathbf{v _j}  \longmapsto \prod _{ j \in \mathbb{Z}} <\mathbf{w}, \mathbf{v _j} >.$$ 
The integer $K$ is called the {\bfi  temporal aggregation period} and the vector $\mathbf{w} $ the {\bfi  temporal aggregation vector}. Notice that the aggregated time series $ Y $ is indexed using the time scale $ \tau = mK$, with $ m \in \mathbb{Z} $ and its components  are given by the $ X ${\bfi  -aggregates} $ X _{ t+K} ^\mathbf{w} $ defined by 
\begin{equation} \label{TAseriesComp}
X _{ t+K} ^\mathbf{w} := w _1 X_{ t+1} +...+ w _K X _{t+K} = Y_{\tau}. 
\end{equation}
\end{definition}
Definition \ref{defTA} can be reformulated in terms of the backward shift operator $ L $ as:
\begin{equation*}
Y = \Pi _K \circ \sum^{K-1}_{i = 0} w _{ K-i} L ^i \left( X \right),
\end{equation*}  
where 
\begin{equation}
\label{projection for ta}
\Pi _K: \enspace \prod_{ j \in \mathbb{Z}} \mathbb{R} \longrightarrow \prod _{ j \in \mathbb{Z}} \mathbb{R}, \enspace
 \left( X_j \right)_{ j \in \mathbb{Z} } \longmapsto \left( X_{K j+1} \right) _{ j \in \mathbb{Z}},
\end{equation}
and the indices of the components $ \left( Z_j \right) _{ j \in \mathbb{Z}} $ of $ Z := \sum^{K-1}_{i = 0} w _{ K-i} L ^i \left( X \right) $ are uniquely determined by the choice $ Z _1 := w _1 X_1 +...+w _K X _K $. 

There are four important particular cases covered by Definition \ref{defTA}, namely:
\begin{description}
\item [(i)] {\bfi  Stock aggregation} (also called systematic sampling, skip-sampling, point-in-time sampling): it is obtained out of \eqref{TAseries} or \eqref{TAseriesComp} by setting $ \mathbf{w} = \left( 0, 0, ...,0, 1 \right)' $.
\item [(ii)] {\bfi  Flow aggregation}: $ \mathbf{w} = \left( 1, 1, ..., 1 \right)' $.
\item [(iii)] {\bfi  Averaging}: $ \mathbf{w} = \left( 1/K, 1/K , ..., 1/K \right)' $.
\item [(iv)] {\bfi  Weighted averaging}: $ \mathbf{w} = \dfrac{1}{K} \left( \xi  _1 , ..., \xi _K \right)' $  such that $ \xi _1 + ...+ \xi _K = 1$. 
\end{description} 
\subsection{Multistep approach to the forecasting of linear temporal aggregates}
Let $ X $ be a time series and $ \mathbf{w} = \left( w _1 , ..., w _K \right)' $ an $ K $-period aggregation vector. Given a finite time realization $\xi _T= \left\{ x_1 , ..., x _T \right\} $ of $ X $ such that $ T=MK $ with $ M \in \mathbb{N} $, we aim at forecasting the aggregate $ w _1 X _{ T+1} + ... + w _K X _{ T+K} $. There are two obvious ways to carry this out; first, we can produce a multistep forecast $ \widehat{ X _{ T+1}}, ..., \widehat{ X _{ T+K}} $ for $ X $ out of which we can obtain the forecast of the aggregate by setting $ \widehat{ X _{ T+K} ^\mathbf{w} }:=w _1 \widehat{ X _{ T+1}}+...+ w _K \widehat{ X _{ T+K}} $. Second, we can temporally aggregate $ X $ using \eqref{TAseries} into the time series $ Y $ given by 
\begin{equation*} 
Y = \mathbb{I} _\mathbf{w} \circ p _K \left( X \right)
\end{equation*}  
and produce a one-step forecast for $ Y $. The following result recalls a well known comparison~\cite{Amemiya1972, Luetkepohl84a, Lutkepohl1986, Lutkepohl1989} of the forecasting performances of the two schemes that we just described using the mean square characteristic error as an optimality criterion. In that setup, given an information set encoded as a $ \sigma $-algebra, the optimal forecast is given by the conditional expectation with respect to it~\cite[page 72]{MR1278033}.   Given a time series $X$, we will denote in what follows by $\sigma \left( \underline{{X_T}} \right) $ the information set generated by a realization $\xi_T=\{ x _1 , \ldots, x _T\} $  of length $T$ of $X$ and the preset $I$ used to produce it; more specifically
\begin{equation*}
\sigma \left( \underline{{X_T}} \right):= \sigma(I\cup \{ x _1 , \ldots, x _T\} ).
\end{equation*}
\begin{proposition}
\label{comparison of frequency predictors}
Let $ X $ be a time series and $ \mathbf{w} = \left( w _1 , ..., w _K \right)' $ a $ K $-period aggregation vector. Let $ Y = \mathbb{I} _\mathbf{w} \circ p _K \left( X \right) $ be the corresponding temporally aggregated time series. Let $ T=MK $ with $ M, T \in \mathbb{N} $ and consider $ \mathcal{F} _T = \sigma \left( \underline{{X_T}} \right) $, $ \mathcal{F} _{ M} = \sigma \left( \underline{Y _{ M}} \right) $  the information sets associated to two histories of $ X $ and $ Y $ of length $ T $ and $ M $, respectively, related to each other by temporal aggregation. Then:
\begin{equation} 
\label{prop31}
{\rm MSFE} \left( E \left[ X _{ T+K} ^\mathbf{w} | \mathcal{F} _{ T} \right] \right) \le {\rm MSFE} \left( E \left[ Y_{ M+1} | \mathcal{F} _{ M} \right] \right).
\end{equation}
\end{proposition}

\begin{remark}
\normalfont
The inequality~(\ref{prop31}) has been studied in detail in the VARMA context by H. L\"utkepohl~\cite{Lutkepohl1986, Lutkepohl1987, Lutkepohl2009} who has fully characterized the situations in which the two predictors are identical and hence have exactly the same performance. This condition is stated and exploited in Section~\ref{Numerical results}, where we illustrate with specific examples the performance of the forecasting scheme that we present in the following pages. 
\end{remark}

In the next two results we spell out the characteristic and the total errors associated to a multistep approach to the forecast of linear aggregates. The characteristic error is given in Proposition~\ref{prop_char_er_aggr} and the total error is provided in Theorem~\ref{theor2} under the same independence hypothesis between the samples used for estimation and forecasting that were already invoked in Theorem~\ref{theor1}.

\begin{proposition}
\label{prop_char_er_aggr}
Let $ X $ be a time series model as in~(\ref{31}), $r=\max\{p,q\}$, $ K $ a temporal aggregation period,  $ \mathbf{w} = \left( w _1, ..., w _K \right)' $ an aggregation vector, and $ \mathcal{F} _T:= \sigma \left(\underline{X _T } \right)$ the information set generated by a realization $\xi_T=\{ x _1 , \ldots, x _T\} $  of length $T$ of $X$. Let $ \widehat{ X _{ T+K} ^\mathbf{w} } $ be the forecast of $ X $-aggregate $ X _{ T+K} ^\mathbf{w} := \sum^{K}_{i = 1} w _i X _{ T+i} $ based on $ \mathcal{F} _T$ using the forecasting scheme in Proposition~\ref{prop1}. Then:
\begin{enumerate}
\item[\bf{(i)}] The forecast $ \widehat{ X _{ T+K} ^\mathbf{w} } $ is given by:
\begin{equation} \label{prop21}
\widehat{ X _{ T+K} ^\mathbf{w}} = \sum^{K}_{i = 1} w _i \sum^{T+i-1+r}_{j = i} \psi _j \varepsilon _{ T+i-j}. 
\end{equation}
\item[\bf{(ii)}] The corresponding mean square forecasting characteristic error is:
\begin{equation} \label{prop22}
{\rm MSFE} \left( \widehat{ X _{ T+K} ^\mathbf{w}} \right)  = E \left[ \left(  X _{ T+K} ^\mathbf{w} - \widehat{ X _{ T+K} ^\mathbf{w}}  \right) ^2 \right] = \sigma ^2 \left[ \sum^{K}_{i = 1} w _i ^2  \sum^{i-1}_{l = 0} \psi _l ^2 + 2 \sum^{K-1}_{i = 1} \sum^{K}_{j = i + 1} w _i w _j  \sum^{i-1}_{l = 0} \psi _l  \psi _{ j-i+l} \right].
\end{equation}
\end{enumerate}
\end{proposition}

\begin{example}{Forecast of stock temporal aggregates.}
\normalfont
It is a particular case of the statement in Proposition~\ref{prop_char_er_aggr} obtained by taking $ \mathbf{w} = \left( 0, ..., 0,1 \right)'$. In this case 
\begin{equation*} 
\widehat{ X _{ T+K} ^\mathbf{w}} = \widehat{ X _{ T+K}} = \sum^{T+K-1+r}_{j= K} \psi _j \varepsilon _{ T+K-j}.
\end{equation*} 
This shows that the forecast of the stock temporal aggregate coincides with the $ K $-multistep forecast of the original time series. Consequently, it is easy to see by using \eqref{prop22} and \eqref{eq:prop13} that 
\begin{equation*} 
{\rm MSFE} \left( \widehat{ X _{ T+K} ^\mathbf{w}} \right)  = {\rm MSFE} \left( \widehat{ X _{ T+K}} \right).
\end{equation*} 
\end{example}
\begin{example}{Forecast of flow temporal aggregates.}
\normalfont
It is a particular case of the statement in Proposition~\ref{prop_char_er_aggr} obtained by taking $ \mathbf{w} = \left( 1, ..., 1 \right)'$. In this case 
\begin{equation*} 
\widehat{ X _{ T+K} ^\mathbf{w}} = \sum^{K}_{i = 1} \sum^{T+i-1+r}_{j= i} \psi _j \varepsilon _{ T+i-j}.
\end{equation*} 
Consequently,
\begin{align}
{\rm MSFE} \left( \widehat{ X _{ T+K} ^\mathbf{w}} \right)  &  = \sigma ^2 \left[ \sum^{K-1}_{j = 0} \left( K-j \right) \psi _j ^2 + 2 \sum^{K-1}_{i = 1} \sum^{K}_{j = i + 1}  \sum^{i-1}_{l = 0} \psi _l  \psi _{ j-i+l} \right] \nonumber \\
&  = \sigma ^2 \left[ \sum^{K-1}_{j = 0} \left( K-j \right) \psi _j ^2 + 2 \sum^{K-1}_{i = 1} \sum^{K-1}_{j = i}  \left( K-j \right)  \psi _{ j-i} \psi _{ j} \right].\nonumber 
\end{align}
\end{example}
\begin{theorem}[Multistep forecasting of linear temporal aggregates]
\label{theor2}
Consider a sample $ \xi_T = \left\{ x_1, \dots , x_T\right\} $ obtained as the realization of a model of the type \eqref{31} using the preset $I$. In  order to forecast out of this sample, we first estimate the parameters of the model $ \widehat{\mathbf{\Psi}} = \left\{ \widehat{\psi} _0 ,  \widehat{\psi}_1, \dots \right\}  $ and $ \widehat{\mathbf{\Pi}} = \left\{ \widehat{\pi} _0 ,  \widehat{\pi}_1, \dots \right\}  $ based on another sample $ \xi'_T$ of the same size that we assume to be independent of $\xi_T$. Let $ \mathbf{w} = \left(w _1 , \dots, w _K \right) '  $ be an aggregation vector and let $\widehat{ X _{ T+K} ^\mathbf{w}}$ be the forecast of the aggregate ${ X_{ T+K} ^\mathbf{w}} := \sum^{K}_{h = 1} w _h X_{T+h}$ based on $\mathcal{F}_T:= \sigma \left(I\cup \xi _T \right) $ using Proposition~\ref{prop_char_er_aggr} and the estimated parameters  $ \widehat{\mathbf{\Psi}}$,  $ \widehat{\mathbf{\Pi}}$. Then: 
\begin{enumerate}
\item[\bf{(i)}]
The optimal forecast $\widehat{X_{T+K}^\mathbf{w}} $ for $X_{T+K}^\mathbf{w}$ given the samples $\xi_T$ and $ \xi'_T$ is
\begin{equation} \label{eq:theor41}
\widehat{ X_{T + K}^\mathbf{w}} = \sum^{K}_{h = 1} w_h \sum^{T+h-1+r}_{j = h} \widehat{ \psi}_{j} {\tilde{\varepsilon}}_{T + h - j},
\end{equation} 
where $ r= \max \{p,q \}$ and $ \tilde{\varepsilon }_t = \sum^{t+r-1}_{j = 0} \widehat{\pi}_j x_{t-j}$.
\item[\bf{(ii)}]
The mean square forecasting error associated to this forecast is
\begin{align} 
\label{eq:theor421}
{\rm MSFE} \left( \widehat{ X_{T + K}^\mathbf{w} }\right) = \sigma ^2 < \mathbf{w}, \left( A + B + C \right) \mathbf{w}>,
\end{align}
where $A$, $B$, $C$ are the matrices with components given by
\begin{equation} 
\label{eq:theor422}
A_{hh'} = \sum^{P(h)}_{i = 0} \sum^{P(h')}_{i' = 0} \psi _i  \psi _{i'} \delta _{h-i, h'-i'},
\end{equation}
\begin{equation} 
\label{eq:theor423}
B_{hh'} = -2 \sum^{P(h)}_{l = 0} \sum^{P(h')}_{i = h'} \sum^{P(h)-i}_{j = 0} \sum^{P(h)-i-j}_{k = 0} \psi _l \psi _{k} E \left[ \widehat{ \psi }_i  \widehat{\pi} _j \right] \delta _{h-l, h'-i-j-k},
\end{equation}
\begin{equation} 
\label{eq:theor424}
C_{hh'} = \sum^{P(h)}_{i = h} \sum^{P(h)-i}_{j = 0} \sum^{P(h)-i-j}_{k = 0} \sum^{P(h')}_{i' = h'} \sum^{P(h')-i'}_{j' = 0} \sum^{P(h')-i'-j'}_{k' = 0} \psi _k \psi _{k'} E \left[ \widehat{ \psi }_i  \widehat{\pi} _j \widehat{\psi}_{i'} \widehat{\pi} _{j'} \right]  \delta _{h-i-j-k,h'-i'-j'-k'},
\end{equation}
with $P(h) = T+h-1+r$, $P(h')=T+h'-1+r$. Notice that $A_{hh'} = A_{hh'}^{char} + A_{hh'}^{res}$,
where 
\begin{equation*}
A_{hh'}^{char}:= \sum^{h-1}_{i = 0}  \sum^{h'-1}_{i' = 0} \psi _i \psi _{i'} \delta _{h-i,h'-i'}, \quad A_{hh'}^{res} := \sum^{P(h)}_{i = h} \sum^{P(h')}_{i' = h'}\psi _i \psi _{i'} \delta _{h-i,h'-i'},
\end{equation*}
and $\sigma ^2 <\mathbf{w},A^{char} \mathbf{w}>$ is the characteristic forecasting error in part {\rm {\bf (ii)}} of  Proposition~\ref{prop_char_er_aggr}.
\item[\bf{(iii)}] Let $\mathbf{\Xi}_{P} := \left( \psi _1 , \dots, \psi_{P(K)} , \pi _1 , \dots, \pi _{P(K)} \right) $, with $P (K)= T+K-1+r $. Using the notation introduced in Lemmas~\ref{delta method} and~\ref{xis and so on} and discarding higher order terms in $1/\sqrt{T} $, the MSFE in \eqref{eq:theor421} can be approximated by  
\begin{equation} \label{eq:theor44}
{\rm MSFE} \left( \widehat{ X_{T + K}^\mathbf{w} }\right) = \sigma ^2 < \mathbf{w}, \left( A^{char} + D + F + G \right) \mathbf{w}>,
\end{equation}
where
\begin{equation*}
A_{hh'}^{char}:= \sum^{h-1}_{i = 0}  \sum^{h'-1}_{i' = 0} \psi _i \psi _{i'} \delta _{h-i,h'-i'},
\end{equation*}
\begin{equation*}
D_{hh'}=\dfrac{1}{T} \sum^{P(h)}_{i = h} \sum^{P(h')}_{i' = h'} \left( \Sigma _ {\mathbf{\Xi}_{P}} \right)_{i,i'}   \delta _{h-i,h'-i'},
\end{equation*}
\begin{equation*}
F_{hh'} = \dfrac{2}{T} \sum^{P(h)}_{i=h} \sum^{P(h')}_{i' = h'} \sum^{P(h')-i'}_{j' = 0} \sum^{P(h')-i'-j'}_{k' = 0} \psi _{i'} \psi _{k'} \left( \Sigma _{{\Xi}_P}\right)_{i, P(K)+j'}  \delta _{h-i, h'-i'-j'-k'},
\end{equation*}
\begin{equation*}\!\!\!\!\!\!\!\!\!\!\!\!\!\!\!\!\!\!G_{hh'} = \dfrac{1}{T}\sum^{P(h)}_{i = h} \sum^{P(h)-i}_{j = 0} \sum^{P(h)-i-j}_{k = 0} \sum^{P(h')}_{i' = h'} \sum^{P(h')-i'}_{j' = 0} \sum^{P(h')-i'-j'}_{k' = 0} \psi _i \psi _{k} \psi_{i'} \psi _{k'}   \left( \Sigma _ {\mathbf{\Xi}_{P} }\right)_{j+P(K),j'+P(K)}    \delta _{h-i-j-k,h'-i'-j'-k'}.
\end{equation*}
\end{enumerate} 

\end{theorem}
\begin{remark}
\normalfont
In order to compute the total error in \eqref{eq:theor44} it is necessary to determine the covariance matrix $\Sigma_{\mathbf{\Xi}_{P} }$. By Lemma~\ref{xis and so on}, it can be obtained  out of the covariance $\Sigma_{\boldsymbol{\beta}} $  matrix associated to the estimator of the ARMA parameters combined with the Jacobian $J_{\Xi_P}$. Details on how to algorithmically compute this Jacobian  are provided in Appendix~\ref{computation of the jacobian}.
\end{remark}

\begin{remark}
\normalfont
Notice that all the matrices involved in the statement of Theorem~\ref{theor2} are symmetric except for $B$ and $F$.
\end{remark} 

\subsection{A hybrid forecasting scheme using aggregated time series models}
\label{A hybrid forecasting scheme using aggregated time series models}

In the previous subsection we presented a forecasting method for linear temporal aggregates based exclusively on the use of high frequency data and models. The performance of this approach has been compared in the literature~\cite{Lutkepohl1986, Lutkepohl1987} with the scheme that consists of using models estimated using the aggregated low frequency data; as it could be expected due to the resulting smaller sample size, this method yields a performance that is strictly inferior to the one based on working in the pure high frequency setup. 

In this section, we introduce and compute the performance of a hybrid recipe that consists of estimating first the model using the high frequency data so that we can take advantage of larger sample sizes and of the possibility of updating the model as new high frequency data become available. This model and the data used to estimate it are subsequently aggregated and used for forecasting. We will refer to this approach as the {\bfi  hybrid forecasting scheme}. The main goal in the following pages is writing  down explicitly the total MSFE associated to this forecasting strategy so that we can compare it using Theorem~\ref{theor2} with the one obtained with the method based exclusively on the use of high frequency data and models. In the next section we use the resulting formulas in order to prove that there are situations in which the hybrid forecasting scheme provides optimal performance for various kinds of temporal aggregation.

The main tool at the time of computing the MSFE associated to the hybrid scheme is again the use of the Delta Method~\cite{Serfling1980} in order to establish the asymptotic normality of the estimation scheme resulting from the combination of high frequency data with the subsequent model temporal aggregation.

In order to make this more explicit, consider a  time series model
$ X $ determined by the parameters $ \boldsymbol{\beta  _X} $ for which an asymptotically normal  estimator $ \widehat{ \boldsymbol{\beta _X} } $ is available, that is, there exists a covariance matrix $ \Sigma _{\boldsymbol{  \beta _X  }}$ such that 
\begin{equation*}
\sqrt{T} \left( \widehat{ \boldsymbol{\beta _X} } - \boldsymbol{\beta _X} \right) \xrightarrow[T \rightarrow \infty]{\enspace dist \enspace}  N \left( 0, \Sigma _{  \boldsymbol{\beta _X}  } \right) 
\end{equation*} 
with $ T $ being the sample size on which the estimation is based. Now, let $ K \in \mathbb{N} $ be an aggregation period, $ \mathbf{w} \in \mathbb{R} ^K $ an aggregation vector, and $ Y := \mathbb{I} _\mathbf{w} \circ p _K \left( X \right) $ the linear temporally aggregated process corresponding to $ X $ and $\mathbf{w}$. 
\begin{proposition}
\label{aggregated covariance matrix}
In the setup that we just described, suppose that the temporally aggregated process $ Y $ is also a parametric time series model and that the parameters $ \boldsymbol{\beta _{ Y}} $ that define it can be expressed as a $ C ^1 $ function $ \boldsymbol{\beta _{ Y}} \left( \boldsymbol{\beta _{ X}} \right) $ of the parameters $ \boldsymbol{\beta _{ X}} $ that determine $ X $. Using the estimator $ \widehat{ \boldsymbol{\beta _{ X}}} $, we can construct an estimator $ \widehat{\boldsymbol{ \beta _{ Y}}} $ for $\boldsymbol{ \beta _{ Y}} $ based on disaggregated $ X $ samples by setting $ \widehat{\boldsymbol{ \beta _{ Y}}}:=\boldsymbol{\beta _{ Y}} \left( \widehat{\boldsymbol{ \beta _{ X}}} \right) $ .
Then 
\begin{equation} \label{prop51}
\sqrt{T} \left( \widehat{ \boldsymbol{\beta _Y} } - \boldsymbol{\beta _Y} \right) \xrightarrow[T \rightarrow \infty]{\enspace dist \enspace}  N \left( 0, \Sigma _{   \boldsymbol{\beta _Y} }  \right), 
\end{equation}
where $ T $ is the disaggregated sample size and 
\begin{equation} \label{prop52}
\Sigma _{ {\boldsymbol{\beta _{ Y}}} } = J _{\boldsymbol{\beta _{ Y}}} \Sigma _{ \boldsymbol{\beta _{ X}} } J _{\boldsymbol{\beta _{ Y}}} ^{ T},
\end{equation}
with $ \left( J _{\boldsymbol{\beta _{ Y}}} \right)  _{ i j} = \dfrac{\partial{ \left(\beta _{ Y} \right) _i}}{\partial{ \left(\beta _{ X}\right)_j }}$ the Jacobian matrix corresponding to the function $ \boldsymbol{\beta _{ Y}} \left( \boldsymbol{\boldsymbol{\beta _{ X}}} \right) $.
\end{proposition}

Once the model temporal aggregation function and its Jacobian have been determined, this proposition can be used to formulate an analog of Theorem~\ref{theor2} for the hybrid forecasting scheme by mimicking the proof of Theorem~\ref{theor1}; the only necessary modification  consists of replacing the asymptotic covariance matrix $\Sigma_{\boldsymbol{\beta _{ X}}} $ of the estimator for the disaggregated model by that of the aggregated model $\Sigma_{\boldsymbol{\beta _{ Y}}} $ obtained using Proposition~\ref{aggregated covariance matrix}. 

We make this statement explicit in the following theorem and then describe how to compute the model aggregation function $\boldsymbol{\beta _{ Y}}(\boldsymbol{\beta _{ X}})$ and its Jacobian $J_{\boldsymbol{\beta _{ Y}}}$ in order to make it fully functional. The construction of these objects is carried out in the ARMA context where the model aggregation question has already been fully studied. Even though all necessary details will be provided later on in the section, all we need to know  at this stage in order to state the theorem is that the linear temporal aggregation of an ARMA(p,q) model is another ${\rm ARMA(p,q }^\ast {\rm )} $ model where
\begin{equation}
\label{new q}
q ^\ast :=\left\lfloor \dfrac{K \left( p+1 \right) +q-p-K^*}{K} \right\rfloor,
\end{equation}
$K$ is the temporal aggregation period, $K ^\ast  $ is the index of the first nonzero entry of the aggregation vector, and the symbol $\lfloor \cdot \rfloor $ denotes the integer part of its argument. We emphasize that if the innovations that drive the disaggregated model are independent with variance $\sigma^2$, this is not necessarily the case for the resulting aggregated model, whose innovations may be only uncorrelated with a different variance $\sigma_\ast ^2$, making it into a so called weak ARMA model.

\begin{theorem}[Hybrid forecasting of linear temporal aggregates]
\label{Hybrid forecasting of linear temporal aggregates}
Let $ \xi_T = \left\{ x_1, \dots , x_T\right\} $ be a sample obtained as a realization of a causal and invertible {\rm ARMA(p,q)}   model $X$ as  in~\eqref{31}.
Let $ \mathbf{w} = \left(w _1 , \dots, w _K \right) '  $ be a temporal aggregation vector such that $T=MK $, for some $M \in \mathbb{N}$, and let $Y = \mathbb{I} _\mathbf{w} \circ p _K \left( X \right)$ be the temporal aggregation of the model $X$ and $ \eta _M:=\{ y _1, \ldots, y _M\} $ the temporal aggregated sample obtained out of $\xi_T $.

We forecast the value of the temporal aggregate $Y_{M+1}= X_{T+K}^\mathbf{w}$ out of the sample $\eta_M $ by first estimating the parameters $\widehat{\boldsymbol{\beta _{ X}} }$ of the model $X$  using  another disaggregated sample $ \xi'_T$ of the same size, that we assume to be independent of $\xi_T$. 
Let $\boldsymbol{\beta _{ Y}}(\boldsymbol{\beta _{ X}})$ be the function that relates the {\rm ARMA} parameter values of the disaggregated and the aggregated model and let  $J_{\boldsymbol{\beta _{ Y}}}$ be its Jacobian. Consider the ${\rm ARMA(p,q }^\ast {\rm )} $ model, with ${\rm q} ^\ast  $ as in~(\ref{new q}), determined by the parameters $\widehat{\boldsymbol{\beta _{ Y}} }:=\boldsymbol{\beta _{ Y}}(\widehat{\boldsymbol{\beta _{ X}} }) $. Then:
\begin{description}
\item [(i)] The optimal forecast $\widehat{Y_{M+1}} $ of the temporal aggregate $Y_{M+1}= X_{T+K}^\mathbf{w}$ based on the information set $\mathcal{F}_M:= \sigma \left(I_K \cup \eta_M \right) $ using Proposition~\ref{prop_char_er_aggr} and the estimated parameters  $\widehat{\boldsymbol{\beta _{ Y}} }$ is given by
\begin{equation}
\label{hybrid forecast}
\widehat{Y_{M+1}}= \sum^{T+r^\ast }_{j = 1} \widehat{ \psi}_{j} {\tilde{\varepsilon}}_{T + h - j},
\end{equation}
where $I _K $ is the preset obtained out of the temporal aggregation of $I$, $ r^\ast = \max \{{\rm p,q} ^\ast  \}$, and $ \tilde{\varepsilon }_t := \sum^{t+r^\ast -1}_{j = 0} \widehat{\pi}_j y_{t-j}$, with $ \widehat{\mathbf{\Psi}} = \left\{ \widehat{\psi} _0 ,  \widehat{\psi}_1, \dots \right\}  $ and $ \widehat{\mathbf{\Pi}} = \left\{ \widehat{\pi} _0 ,  \widehat{\pi}_1, \dots \right\}  $ the parameters corresponding to the causal and invertible representations of the temporally aggregated {\rm ARMA}  model with parameters $\widehat{\boldsymbol{\beta _{ Y}} }$.
\item[\bf{(ii)}]
Let $\mathbf{\Xi}_P := \left( \psi _1 , \dots, \psi_P , \pi _1 , \dots, \pi _P \right)' $ with $P = T+r^\ast  $.   Discarding higher order terms in $1/\sqrt{T}$, the MSFE corresponding to the forecast~\eqref{hybrid forecast} can be approximated by  
\begin{align} 
\label{hybrid error}
\!\!\! {\rm MSFE} \left( \widehat{ Y_{M+1} }\right) &= \sigma ^2_\ast   + \sigma ^2 _\ast\dfrac{1}{T}\Biggl[  \sum^{P}_{i = h}(\Sigma _{\mathbf{\Xi}_P})_{i,i}  + 2 \sum^{P}_{i = h} \sum^{P-i}_{j = 0} \sum^{P-i-j}_{k = 0}\psi_{i} \psi_k  (\Sigma _{\mathbf{\Xi}_P})_{i+j+k,j+P} \nonumber \\
& + \sum^{P}_{i = 1} \sum^{P-i}_{j = 0} \sum^{P-i-j}_{k = 0} \sum^{P}_{i' = 1} \sum^{P-i'}_{j' = 0} \sum^{P-i'-j'}_{k' = 0} \psi _i \psi _k \psi _{i'} \psi _{k'}  (\Sigma _{\mathbf{\Xi}_P})_{j+P,j'+P} \delta _{i+j+k,i'+j'+k'} \Biggl],
\end{align} 
where $ \sigma ^2_\ast $ is the variance of the innovations of the aggregated ${\rm ARMA(p,q }^\ast {\rm )} $ model $Y$ and $\Sigma_{\boldsymbol{\Xi}_P} $ is the covariance matrix given by
\begin{equation*}
\Sigma_{\boldsymbol{\Xi}_P} =J _{\mathbf{\Xi}_P} \Sigma _{\boldsymbol{\beta_Y}} J_{\mathbf{\Xi}_P}'=J _{\mathbf{\Xi}_P} J _{\boldsymbol{\beta _{ Y}}} \Sigma _{ \boldsymbol{\beta _{ X}} } J _{\boldsymbol{\beta _{ Y}}} ^{ T} J_{\mathbf{\Xi}_P}',
\end{equation*}
with $J _{\boldsymbol{\beta _{ Y}}} $  the Jacobian matrix corresponding to the function $ \boldsymbol{\beta _{ Y}} \left( \boldsymbol{\beta _{ X}} \right) $ and $ J_{\mathbf{\Xi}_P}$ the Jacobian of $\mathbf{\Xi}_P(\boldsymbol{\beta}) := \left( \psi _1(\boldsymbol{\beta_Y}) , \dots, \psi_P(\boldsymbol{\beta_Y}) , \pi _1 (\boldsymbol{\beta_Y}), \dots, \pi _P(\boldsymbol{\beta_Y})\right)' $.
\end{description}
\end{theorem}

As we announced above, we  conclude this section by describing in detail the  parameters aggregation function $\boldsymbol{\beta _{ Y}}(\boldsymbol{\beta _{ X}})$ and its Jacobian $J_{\boldsymbol{\beta _{ Y}}}$, so that all the ingredients necessary to apply formula~(\ref{hybrid error}) are available to the reader. In order to provide explicit expressions regarding these two elements, we provide a brief review containing strictly the  results on the temporal aggregation of ARMA processes that are necessary for our discussion; for more ample discussions about this topic we refer the reader to~\cite{Amemiya1972, Tiao1972, Brewer1973, TiaoWei1976, WeiTempAggr1979, Stram1986, Wei:book, Silvestrini2008} and references therein.

\medskip

\noindent {\bf The temporal aggregation function $\boldsymbol{\beta _{ Y}}(\boldsymbol{\beta _{ X}})$.}
Consider the ARMA(p,q) model $ \mathbf{\Phi} \left( L \right) X = \mathbf{\Theta} \left( L \right) \varepsilon $, where $ \mathbf{\Phi} \left( L \right) = 1 - \sum^{p}_{i = 1} \phi _i L ^i $ and $ \mathbf{\Theta} \left( L \right) = 1+ \sum^{q}_{i = 1} \theta _i L ^i $. In order to simplify the discussion and to avoid hidden periodicity phenomena, we place ourselves in a generic situation in which the autoregressive and moving average polynomials of the model that we want to temporally aggregate have no common roots and all roots are different (see~\cite{Wei:book} for the general case).
Consider now a $ K $-period aggregation vector $ {\bf w} = \left( w _1, ..., w _K \right) ' $. Our first goal is to find polynomials $ \mathbf{T} \left( z \right) $ and $ \mathbf{\Phi^*} \left( z \right) $ that satisfy 
\begin{equation} \label{TLPHIL}
\mathbf{T} \left( L \right) \mathbf{\Phi} \left( L \right) = \mathbf{\Phi^*} \left( L ^K \right) \circ \Pi _K \circ \sum^{K}_{i = 1} w _i L ^i, 
\end{equation} 
with $\Pi _K: \enspace \prod_{ j \in \mathbb{Z}} \mathbb{R} \longrightarrow \prod _{ j \in \mathbb{Z}} \mathbb{R}$ as in~(\ref{projection for ta}). The intuition behind  \eqref{TLPHIL} is that for any time series $ X $, its temporal aggregation $ Y = \Pi _K \circ \sum^{K}_{i = 1} w _i L ^i \left( X \right) $ satisfies 
\begin{equation} \label{TLPHIL1}
\mathbf{T} \left(  L \right) \mathbf{\Phi} \left( L \right) X = \mathbf{\Phi^*} \left( L ^K \right) Y.
\end{equation}
In other words, the polynomial $ \mathbf{T} \left( L \right) $, that we will call the {\bfi   temporal aggregation polynomial} transforms the AR polynomial for $ X $ into an AR polynomial for $ Y $ in the aggregated time scale units. Let 
\begin{equation*}
\mathbf{T} \left( L \right) = t _0 + t _1 L + ... + t _n L ^n \quad \mbox{and} \quad \mathbf{\Phi^*} \left( L ^K \right) = 1 - \phi _1 ^* L ^K - ... - \phi _c ^*  L ^{ Kc} 
\end{equation*}
 be the unknown polynomials. Equation \eqref{TLPHIL} can be written in  matrix form as~\cite{Brewer1973}:
\begin{equation} 
\label{Brewer}
\underbrace{
\left(
\begin{array}{ccccc} 
t_0 & 0 & 0 & \ldots & 0\\ 
t_1 & t_0 & 0 & \ldots & 0\\ 
t_2 & t_1 & t_0 & \ldots & 0\\ 
\vdots & \vdots & \vdots & \ddots & \vdots \\
t_p & t_{p-1} & t_{p-2} & \ldots & t_0 \\
\vdots & \vdots & \vdots & \ddots & \vdots \\ 
t_n & t_{n-1} & t_{n-2} & \ldots & t_{n-p} \\
0 & t_{n} & t_{n-1} & \ldots & t_{n-p-1} \\
0 & 0 & t_{n} & \ldots & t_{n-p-2} \\
\vdots & \vdots & \vdots & \ddots & \vdots \\ 
0 & 0 & 0 & \ldots & t_{n} \\
\end{array} \right)
}_{{\Huge{\rm A}}}
\underbrace{
\left(
\begin{array}{c} 
1\\ 
- \phi _1 \\ 
\vdots \\
- \phi _p \\
\end{array} \right)
}_\text{Z}
=
\underbrace{
\left(
\begin{array}{c} 
\widetilde{ \mathbf{w}}\\ 
- \phi _1 ^* \widetilde{ \mathbf{w}}\\ 
\vdots \\
 \vdots \\ 
\vdots \\
\vdots\\
\vdots\\
\vdots\\
- \phi _c ^* \widetilde{ \mathbf{w}} \\
\end{array} \right)
}_\text{D}
\end{equation}
where $ \widetilde{ \mathbf{w}} = \left( w _K , w _{ K-1}, ..., w _1 \right) '$ is the reflection of $\mathbf{w}$. We start by determining the unknown values $ n $ and $ c $ using the two following dimensional restrictions:
\begin{itemize}
\item Since $A$ is a matrix of size $ \left( n + p+1 \right) \times  \left( p +1 \right) $,  $Z$ and $D$ are vectors of size $p+1$ and $cK +K $, respectively, and we have $ AZ=D$ then, necessarily
\begin{equation} \label{2}
n + p + 1 = cK + K.
\end{equation}
\item 
The system $ AZ=D$ contains $ n+p+1 $ equations that need to coincide with the number of unknowns, that is, the $ n+1+c $ coefficients $(t _0, t _1, \ldots, t _n)$ and $\left( \phi _1 ^*, \ldots, \phi _c ^* \right)  $ of the polynomials  $ \mathbf{T} \left( z \right) $ and $ \mathbf{\Phi^*}(z)$, respectively. Consequently,
\begin{equation} \label{3}
n+p+1=n+c+1.
\end{equation}
\end{itemize}
The conditions \eqref{2} and \eqref{3} yield
\begin{equation} \label{4}
\left\{
\begin{array}{rl} 
c&= p, \\ 
n&= pK+K-p-1 = \left( p+1 \right) \left( K - 1 \right).
\end{array} \right.
\end{equation}
The first condition in \eqref{4} shows that {\it the autoregressive order does not change under temporal aggregation}. 

Let now $ K ^* \le K $ be the index of the first nonzero component in the vector ${\bf  w} $. This implies that $\widetilde{{\bf w}} $ is a vector of the form $ \widetilde{{\bf w}} = \left( w _K , w _{ K-1}, ..., w _{K ^*}, 0 , ..., 0 \right) ' $ with $ w _K, w_{K-1}, \ldots, w _{K ^*} \ne 0 $. Since \eqref{Brewer} is a matrix representation of the polynomial equality in \eqref{TLPHIL}, we hence have that $ \deg \left(\mathbf{T} \left( L \right) \mathbf{\Phi} \left( L \right) \right)= \deg \left(\mathbf{D} \left( L \right) \right)$, where $ \mathbf{D} \left( L \right) $ is the polynomial associated to the vector in the right hand side of \eqref{Brewer}. It is clear that $ \deg \left(\mathbf{D} \left( L \right)\right) = K \left( p+1 \right) - K ^* $; as $ \deg \left(\mathbf{T} \left( L \right) \mathbf{\Phi} \left( L \right) \right)= n+p$,  the degree $ n $ of $ \mathbf{T} \left( L \right) $ is therefore
\begin{equation}
\label{degree of t}
n = K \left( p+1 \right) - p - K ^*.
\end{equation}

Solving the polynomial equalities  \eqref{Brewer}, we have found polynomials $ \mathbf{T} \left( L \right) $ and $ \mathbf{\Phi^*}\left( L ^K \right) $ such that the temporally aggregated time series $ Y $ satisfies 
\begin{equation} \label{5}
\mathbf{\Phi^*} \left( L ^K \right) Y = \mathbf{T} \left( L \right) \mathbf{\Phi} \left( L \right) X 
\end{equation}
Our goal now is showing the existence of a polynomial $ \mathbf{\Theta} ^* \left( L ^K \right) $ and a white noise $ \left\{ \varepsilon ^* \right\} \sim \text{WN} \left( 0, \sigma _{ \ast} ^2 \right) $ such that 
\begin{equation}
\mathbf{T} \left( L \right) \mathbf{\Theta} \left( L \right) \varepsilon _{ lK} = \mathbf{\Theta^*} \left( L ^K \right) \varepsilon _{ lK}^* \enspace \text{for any} \enspace l \in \mathbb{Z}.
\end{equation}
This equality, together with \eqref{5} shows that the temporally aggregated process $ Y $ out of the ARMA process $ X $ is a weak ARMA process, as it satisfies the relation
\begin{equation} \label{6}
\mathbf{\Phi^*} \left( B \right) Y = \mathbf{\Theta^*} \left( B \right) \varepsilon ^*, \enspace \left\{ \varepsilon ^* \right\} \sim \text{WN} \left( 0, \sigma _{ \ast} ^2 \right), 
\end{equation}
where $ B = L ^K $, $ \mathbf{\Phi^*} \left( B \right) $ is a polynomial of degree $p$, and $ \mathbf{\Theta^*} \left( B \right) $ a polynomial of degree $ \left\lfloor \dfrac{n+q}{K} \right\rfloor$ whose coefficients will be determined in the following paragraphs. We recall that the symbol $\lfloor \cdot \rfloor $ denotes the integer part of its argument. Indeed, by \eqref{degree of t} we have that $ \deg \left(\mathbf{T} \left( L \right) \mathbf{\Theta} \left( L \right)\right) = K \left( p+1 \right) -p-K^* +q= n+q$. Additionally, by \eqref{5}
\begin{equation} \label{7}
\mathbf{T} \left( L \right) \mathbf{\Theta} \left( L \right) \varepsilon = \mathbf{T} \left( L \right) \mathbf{\Phi} \left( L \right) X = \mathbf{\Phi^*} \left( L ^K \right) Y.
\end{equation}
Since the right hand side of this expression only involves time steps that are integer multiples of $ K $, the relation~(\ref{7}) only imposes requirements on the left hand side at those time steps. Moreover, it is easy to see that
\begin{equation}
E \left[ \left(\mathbf{T} \left( L \right) \mathbf{\Theta} \left( L \right) \varepsilon _{ lK}\right) \left( \mathbf{T} \left( L \right) \mathbf{\Theta} \left( L \right) \varepsilon _{ lK+jK}\right) \right] = 0,
\end{equation}
for any $ Kj > K \left( p+1 \right) +q-p-K^* $. This implies that  the process is $ \left\{ \mathbf{T} \left( L \right) \mathbf{\Theta} \left( L \right) \varepsilon _{ lK} \right\} _{ l \in \mathbb{Z}} $ is $ \left( K \left( p+1 \right) +q-p-K^*\right)$-correlated, which guarantees in turn by~\cite[Section 3.2]{BrocDavisYellowBook}  the existence of a weak MA($ q^* $) representation
\begin{equation} \label{7}
\mathbf{T} \left( L \right) \mathbf{\Theta} \left( L \right) \varepsilon _{lK} =\mathbf{\Theta^*} \left( B \right) \varepsilon _{ lK} ^*, \enspace l \in \mathbb{Z},
\end{equation}
where 
\begin{equation}
\deg \left(\mathbf{\Theta^*}\left( B \right)\right) = \left\lfloor \dfrac{K \left( p+1 \right) +q-p-K^*}{K} \right\rfloor :=q^*.
\end{equation}
The coefficients of the polynomial $ \mathbf{\Theta^*} \left( B \right) $ are obtained by equating the autocovariance functions of the processes on both sides of \eqref{7} at lags $ 0, K , 2K, ..., q^*K $, which provide $ q^* + 1 $ nonlinear equations that determine uniquely the $ q^* +1 $ unknowns corresponding to the coefficients $ \theta _1 ^*,\ldots, \theta _{ q^*} ^*$ of  the polynomial $\mathbf{\Theta^*} \left( B \right) $ and the variance $\sigma _{ \ast}^2 $ of the white noise of the aggregated model. 

In order to explicitly write down the equations that we just described, let us denote $ \mathbf{C} \left( L \right) := \mathbf{T} \left( L \right) \mathbf{\Theta} \left( L \right) $ as in \eqref{7} and set $ \mathbf{C} \left( L \right) = \sum^{n+q}_{i = 0} c _i L ^i $. Let now $ \gamma $ and $\Gamma $ be the autocovariance functions of the ${\rm MA \left( n+q \right)} $ and ${\rm MA \left( q^* \right) }$ processes
\begin{equation}  \label{8}
V _t = \mathbf{C} \left( L \right) \varepsilon  _t  \quad \mbox{and} \quad U _\tau = \mathbf{\Theta^*} \left( B \right) \varepsilon _{ \tau } ^*, \quad \mbox{respectively},
\end{equation} 
which are given by 
\begin{equation} \label{9}
\gamma \left( i \right) = \sigma _{ \varepsilon } ^2 \sum^{n+q-i}_{l = 0} c _l c _{l+i}\quad \mbox{and} \quad \Gamma \left( j \right) = \sigma _{ \ast} ^2 \sum^{q ^* - j}_{l = 0} \theta _l ^* \theta _{l+j} ^*.
\end{equation}
Consequently, the coefficients of the polynomial $ \mathbf{\Theta^*} \left( B \right) $ and the variance $ \sigma _{ \ast} ^2 $ are uniquely determined by the $ q^* +1$ equations 
\begin{equation} \label{10}
\sigma _{ \varepsilon } ^2 \sum^{n+q-i}_{l = 0} c _l c _{l+jK} = \sigma _{ \ast} ^2 \sum^{q ^* - j}_{l = 0} \theta _l ^* \theta _{l+j} ^*, \enspace j = 0, 1, \dots, q^*. 
\end{equation}
The equations \eqref{10} can be written down in matrix form, which is convenient later on at the time of spelling out the Jacobian of the aggregation function. Indeed, we can write:
\begin{equation} \label{11}
\gamma \left( i \right) = \sigma _{ \varepsilon } ^2 \left[  \overline{\mathbf{T} \left( L \right) \mathbf{\Theta} \left( L \right)} ' S _i   \overline{\mathbf{T} \left( L \right) \mathbf{\Theta} \left( L \right)}  \right] = \sigma _{ \varepsilon } ^2 \left[  \overline{\mathbf{C} \left( L \right) }' S _i  \overline{\mathbf{C} \left( L \right) }  \right] , 
\end{equation}
\begin{equation} \label{11}
\Gamma \left( i \right) = \sigma _{ \ast} ^2 \left[  \overline{ \mathbf{\Theta^*}  \left( B \right) } ' S _i   \overline{ \mathbf{\Theta^*}  \left( B \right) }  \right] ,
\end{equation}
where the bars over the polynomials in the previous expressions denote the corresponding coefficient vectors, that is, given a polynomial $ q \left( x \right) = \sum^{n}_{i = 1} a _i x ^i $, then $ \overline{q \left( x \right) } = \left( a _0, a _1 , \dots, a _n \right)'$. Additionally $ S _i $ is the lower $i$th-shift matrix, that is,
$$
\left( S _i \right) _{jl} = \delta _{j-l, i}.
$$ 
For any given vector $ \mathbf{v}= \left( v _1 , \dots, v _n \right) '$,  $ S _i \mathbf{v}= ( \underbrace{0, \dots, 0}_i, v _1 , \dots, v _{n-i}) '$. With this notation, the equations  \eqref{10} can be rewritten as 
\begin{equation} \label{12}
\sigma _{ \varepsilon } ^2 \left[  \overline{\mathbf{T} \left( L \right) \mathbf{\Theta} \left( L \right)} ' S _{jK}   \overline{\mathbf{T} \left( L \right) \mathbf{\Theta} \left( L \right)}  \right] = \sigma _{ \ast} ^2 \left[  \overline{ \mathbf{\Theta^*}  \left( B \right) } ' S _j   \overline{ \mathbf{\Theta^*}  \left( B \right) }  \right], \enspace j = 0, 1, \dots, q^*. 
\end{equation}
In conclusion,  if we denote $ \boldsymbol{\beta _{ X}} = \left( \mathbf{\Phi} , \mathbf{\Theta} \right)  $ and $\boldsymbol{\beta _{ Y}} = \left( \mathbf{\Phi^*} ,\mathbf{\Theta^*} \right)   $, the construction that we just examined shows that 
\begin{equation}
\label{beta function 1}
\boldsymbol{\beta _{ Y}} \left( \boldsymbol{\beta _{ X}}\right) = \left( \mathbf{\Phi^*} \left( \mathbf{\Phi} \right) ,\mathbf{\Theta^*} \left( \mathbf{\Phi} , \mathbf{\Theta} \right) \right). 
\end{equation}
The function $\mathbf{\Phi^*} \left( \mathbf{\Phi} \right) $ is given by the solution of the polynomial equalities \eqref{Brewer} and $\mathbf{\Theta^*} \left( \mathbf{\Phi} , \mathbf{\Theta} \right) $ by the coefficients $(t _0, t _1, \ldots, t _n)$ determined by \eqref{Brewer} and the solutions of~(\ref{12}).

\begin{example}{Stock temporal aggregation of an {\rm ARMA(p,q)} model.} 
\normalfont
In this case, $\mathbf{w}=(0, \ldots,0,1)'$ and hence $ K^* = K$, $ n = p \left( K-1 \right) $, and $ q ^* = \left\lfloor \dfrac{p \left( K-1 \right) +q}{K}\right\rfloor $.  
\end{example}
\begin{example}{Flow temporal aggregation of an {\rm ARMA(p,q)} model.}
\normalfont
In this case, $\mathbf{w}=(1, \ldots,1)'$ and hence $ K^* = 1$, $ n = \left( p+1 \right)  \left( K-1 \right) $, $ q ^* = \left\lfloor \dfrac{\left( p+1 \right)  \left( K-1 \right) +q}{K}\right\rfloor $.  
\end{example}

\medskip

\noindent {\bf The Jacobian $J_{\boldsymbol{\beta _{ Y}}}$ of the temporal aggregation function $\boldsymbol{\beta _{ Y}} \left( \boldsymbol{\beta _{ X}}\right)  $.}
The goal of the following paragraphs is the computation of the Jacobian $J_{\boldsymbol{\beta _{ Y}}}$ of the function $\boldsymbol{\beta _{ Y}} \left( \boldsymbol{\beta _{ X}}\right) = \left( \mathbf{\Phi^*} \left( \mathbf{\Phi} \right) , \mathbf{\Theta^*} \left( \mathbf{\Phi} , \mathbf{\Theta} \right) \right)  $ in~(\ref{beta function 1}).
We first compute the Jacobian  of the function $ \mathbf{\Phi^*} \left( \mathbf{\Phi} \right)$ by taking derivatives with respect to the components of the vector $ \mathbf{\Phi} $ on both sides of the equations \eqref{Brewer} that determine  $\mathbf{\Phi^*} \left( \mathbf{\Phi} \right) $. This  results in the following $p$  matrix equations
\begin{equation} \label{Brewer_der}
\left(
\begin{array}{ccccc} 
\dfrac{\partial{t_0}}{\partial{\phi _i }} & 0 &  \ldots & 0\\ 
\dfrac{\partial{t_1}}{\partial{\phi _i }} & \dfrac{\partial{t_0}}{\partial{\phi _i }}& \ldots & 0\\ 
\dfrac{\partial{t_2}}{\partial{\phi _i }} & \dfrac{\partial{t_1}}{\partial{\phi _i }} & \ldots & 0\\ 
\vdots & \vdots & \ddots & \vdots \\
\dfrac{\partial{t_p}}{\partial{\phi _i }} &\dfrac{\partial{t_{p-1}}}{\partial{\phi _i }}& \ldots & \dfrac{\partial{t_0}}{\partial{\phi _i }} \\
\vdots & \vdots & \ddots & \vdots \\ 
\dfrac{\partial{t_n}}{\partial{\phi _i }} & \dfrac{\partial{t_{n-1}}}{\partial{\phi _i }} & \ldots & \dfrac{\partial{t_{n-p}}}{\partial{\phi _i }} \\
0 & \dfrac{\partial{t_n}}{\partial{\phi _i }} &  \ldots & \dfrac{\partial{t_{n-p-1}}}{\partial{\phi _i }}\\
0 & 0 & \ldots & \dfrac{\partial{t_{n-2}}}{\partial{\phi _i }} \\
\vdots & \vdots &\ddots & \vdots \\ 
0 & 0 & \ldots&\dfrac{\partial{t_n}}{\partial{\phi _i }} \\
\end{array} \right)
\left(
\begin{array}{c} 
1\\ 
- \phi _1 \\ 
- \phi _2 \\ 
\vdots \\
- \phi _p \\
\end{array} \right)
+\left(\begin{array}{ccccc} 
t_0 & 0 &  \ldots & 0\\ 
t_1 & t_0 &  \ldots & 0\\ 
t_2 & t_1 &  \ldots & 0\\ 
\vdots & \vdots & \ddots & \vdots \\
t_p & t_{p-1} & \ldots & t_0 \\
\vdots & \vdots & \ddots & \vdots \\ 
t_n & t_{n-1} &  \ldots & t_{n-p} \\
0 & t_{n} & \ldots & t_{n-p-1} \\
0 & 0 & \ldots & t_{n-p-2} \\
\vdots & \vdots & \ddots & \vdots \\ 
0 & 0 & \ldots & t_{n} \\
\end{array} \right)
\left(
\begin{array}{c} 
0\\ 
- \dfrac{\partial{\phi _1}}{\partial{\phi _i }} \\ 
- \dfrac{\partial{\phi _2}}{\partial{\phi _i }} \\ 
\vdots \\
- \dfrac{\partial{\phi _p}}{\partial{\phi _i }}  \\
\end{array} \right)
=\left(
\begin{array}{c} 
0\\
- \tilde{w} \dfrac{\partial{\phi _1 ^*}}{\partial{\phi _i }} \\ 
- \tilde{w} \dfrac{\partial{\phi _2 ^*}}{\partial{\phi _i }} \\ 
\vdots \\
- \tilde{w} \dfrac{\partial{\phi _p ^*}}{\partial{\phi _i }} \\
\end{array} \right),
\end{equation}
with $ i = 1, ..., p $. These equations uniquely determine the $(1,1)$-block $ \dfrac{\partial{\mathbf{\Phi^*}}}{\partial{\mathbf{\Phi} }} = \left( \dfrac{\partial{\phi _i ^*}}{\partial{\phi _j }} \right) _{ i,j} $ of the Jacobian $J_{\boldsymbol{\beta _{ Y}}}$, as well as the derivatives 
\begin{equation} \label{derd}
 d_{ ij}:= \dfrac{\partial{t _i }}{\partial{\phi _j}}
 \end{equation}
  that will be needed later on in the computation of the remaining blocks of the Jacobian. 
Given that there is no $ \mathbf{\Theta} $ dependence on the function $\mathbf{\Phi^*} \left( \mathbf{\Phi} \right)$, the $(1,2)$-block of the Jacobian is a zero matrix of size $ p \times  q^* $.

The remaining two blocks are computed by using the function $ \mathbf{\Theta} ^\ast \left( \mathbf{\Phi} , \mathbf{\Theta} \right) $ uniquely determined by the equations \eqref{12}. Its derivatives are obtained out of a new set of equations resulting from the differentiation of both sides of this relation, namely,

\begin{align} \label{13}
 &  \sigma _{\varepsilon } ^2 \left[ \overline{\dfrac{\partial{\mathbf{T} \left( L \right) }}{\partial{\left( \beta _X \right) _i  }} \mathbf{\Theta} \left( L \right) } + \overline{ \mathbf{T} \left( L \right) \dfrac{\partial{\mathbf{\Theta} \left( L \right) }}{\partial{\left( \beta _X \right) _i  }} }\right] ' S_{jK}  \overline{\mathbf{T} \left( L \right) \mathbf{\Theta} \left( L \right) } + \overline{ \mathbf{T} \left( L \right) \mathbf{\Theta} \left( L \right)  } ' S_{jK} \left[ \overline{\dfrac{\partial{\mathbf{T} \left( L \right) }}{\partial{\left( \beta _X \right) _i  }} \mathbf{\Theta} \left( L \right) } + \overline{ \mathbf{T} \left( L \right) \dfrac{\partial{\mathbf{\Theta} \left( L \right) }}{\partial{\left( \beta _X \right) _i  }} }\right]   \nonumber \\
& = \dfrac{\partial{\sigma _{\ast} ^2}}{\partial{\left( \beta _X \right) _i  }}  \left[ \overline{ \mathbf{\Theta^*} \left( B \right)}' S _j \overline{ \mathbf{\Theta^*} \left( B \right)} \right] + 
 \sigma _{\ast} ^2 \left[ \overline{\dfrac{\partial{\mathbf{\Theta^*} \left( B \right) }}{\partial{\left( \beta _X \right) _i  }} } S_{jK}  \overline{\mathbf{\Theta^*} \left( B \right) } + \overline{ \mathbf{\Theta^*} \left( B \right) }' S_{j} \overline{\dfrac{\partial{\mathbf{\Theta^*} \left( B \right) }}{\partial{\left( \beta _X \right) _i  }}} \right] ,\\
& j = 0, 1, \dots, q^*,\enspace   i = 1, \dots, p+q. \nonumber
\end{align}

We recall that the entries of the vector $ \overline{\dfrac{\partial{\mathbf{T} \left( L \right) }}{\partial{\phi} _i }}$ correspond to the values previously obtained in \eqref{derd} and  that $ \overline{\dfrac{\partial{\mathbf{T} \left( L \right) }}{\partial{ \theta} _i }} = \overline{0}$.
Expression \eqref{13} provides $ \left( q^* + 1 \right) \left( p+q \right) $ equations that allow us to find the values of the $ \left( q^* + 1 \right) \left( p+q \right) $ unknowns
\begin{equation}
\dfrac{\partial{\left( \overline{\mathbf{\Theta^*} \left( B \right) }\right) _j  }}{\partial{\phi} _r },  \enspace \dfrac{\partial{\left( \overline{\mathbf{\Theta^*} \left( B \right) }\right) _j  }}{\partial{\theta} _s },  \enspace \dfrac{\partial{\sigma _{\ast} ^2 }}{\partial{\phi} _r }, \enspace \dfrac{\partial{\sigma _{\ast} ^2 }}{\partial{\theta} _s}, \enspace j = 1, \dots, q^\ast , \enspace r = 1, \dots, p, \enspace s = 1, \dots, q.
\end{equation}

\section{Comparison of forecasting efficiencies. Examples.}
\label{Numerical results}
In the previous section we proposed a new hybrid scheme for the forecasting of temporal aggregates coming from ARMA processes. We recall that this strategy consists of first using high frequency disaggregated data for estimating a model; then we temporally aggregate both the data and the model, and finally we forecast using these two ingredients. As we announced  in the introduction, there are situations in which our strategy is optimal with respect to the total error, that is, the predictor constructed  following this procedure performs better than the one based exclusively on high frequency  data and the underlying disaggregated model. In this section we give a  few examples of specific models for which our scheme provides optimal efficiency of prediction.  Before we proceed, we introduce abbreviations for the various predictors that we will be working with:
\begin{description}
\item[(i)] \ \,{\bf Temporally aggregated multistep predictor (TMS predictor)}: this is the denomination that we use for the forecast of the aggregate that is constructed out of the disaggregated data and the underlying disaggregated model estimated on them.
\item[(ii)] \ {\bf Temporally aggregated predictor (TA predictor)}: this is the forecast based on use of the temporally aggregated sample and a model estimated on it.
\item[(iii)] {\bf Hybrid predictor (H predictor)}: this is the predictor introduced in Section~\ref{A hybrid forecasting scheme using aggregated time series models} whose performance is spelled out in Theorem~\ref{Hybrid forecasting of linear temporal aggregates}. In this scheme, a first model is estimated on the disaggregated high frequency data sample, then the data and the model are temporally aggregated with an  aggregation period that coincides with the forecasting horizon; finally, both the temporally aggregated model and the sample are used to produce a one-step  ahead forecast that amounts to a prediction of the aggregate we are interested in. 
\item[(iv)]
{\bf  Optimal hybrid predictor (OH predictor)}: this predictor is constructed by taking the  multistep implementation of the H predictor that yields the smallest total error. More explicitly, suppose that the aggregate that we want to forecast involves $K$ time steps; let $\{K _1, \ldots, K _r\} $ be the positive divisors of $K$ and $\{C _1, \ldots, C _r\} $ the corresponding quotients, that is, $K=K _i C _i $ for each $i \in \{1, \ldots, r\} $. There are aggregation schemes (stock and flow for example) for which a $K$-temporal aggregate can be obtained as the aggregation of $C _i$ $K _i $-temporal aggregates, for all $i \in \{1, \ldots, r\} $. The total error associated to the forecasting of these aggregates using a multistep version of the H predictor obviously depends on the factor $K _i $ used. The OH predictor is the one associated to the factor $K _i $ that minimizes the total error.
\end{description}

As we already mentioned, the forecasting performance of the TMS predictor is always superior or equal than that of the TA predictor when we take into account only the characteristic error, and it is strictly superior when the total error is considered. In view of these results and given that the H and the OH predictors carry out the forecasting with temporally aggregated data, they are going to produce worse characteristic errors than their TMS counterpart; hence, the only way in which the H and OH predictors can be competitive in terms of total error performance is by sufficiently lowering the estimation error. In order to check that they indeed do so, we will place ourselves in situations that are particularly advantageous in this respect and will choose models for which the TMS and the TA predictors have identical characteristic errors and hence it is only the estimation error that makes a difference as to the total error. The linear models for which this coincidence of characteristic errors takes place have been identified in the works of H. L\"utkepohl~\cite{Lutkepohl1986, Lutkepohl1987, Lutkepohl2009} via the following statement.

\begin{theorem}[L\"utkepohl]
Let $X_t = \sum^{\infty}_{i = 0} \psi_i  \varepsilon_{t - i}$ be a linear causal process and let $\mathbf{w}=\left( w_1, \dots, w_K \right)' \in \mathbb{R}^K$ be a $K$-period aggregation vector. Then the TMS and TA predictors for the $K$-temporal aggregate determined by $\mathbf{w} $ have identical associated characteristic errors if and only if the following identity holds:
\begin{equation} \label{charerrorequal}
\left( \sum^{K-1}_{i = 0} w_{K-i} L^i \right) \mathbf{\Psi} \left( L \right) = \left( \sum^{\infty}_{j = 0} \left( \sum^{K-1}_{i = 0} w_{K-i}  \psi_{jK-i} \right) L^{jK} \right) \left( \sum^{K-1}_{j = 0} \left( \sum^{j}_{i = 0} w_{K-i}  \psi_{j-i} \right) L^{j} \right).
\end{equation} 
\end{theorem}

The equality \eqref{charerrorequal} is satisfied for both stock ($\mathbf{w}=\left( 0, \dots, 0, 1 \right)'$) and flow aggregation ($\mathbf{w}=\left( 1, \dots, 1 \right)'$) if $\{X_t\}$ is a purely seasonal process with period $K$, that is,
\begin{equation}
\label{periodic condition}
X_t = \sum^{\infty}_{i = 0} \psi_{iK}\varepsilon_{t-iK}. 
\end{equation}

Given a specific model we want to compare the performances of the H and the TMS predictors for a variety of forecasting horizons. Given that condition \eqref{charerrorequal} is different for each aggregation period $K$ and cannot be solved simultaneously for several of them, we will content ourselves either with approximate solutions that are likely to produce very close H and TMS characteristic errors for several periods $K$ or with exact solutions that provide exactly equal errors for only a prescribed aggregation period. The following points describe how we have constructed examples following the lines that we just indicated:
\begin{itemize}
\item
We first choose the orders $p$ and $q$ of the disaggregated ARMA(p,q) model that we want to use as the basis for the example.
\item
We fix an aggregation period $K$ and a number $n$ of  parameters $\psi_i$ for which the equation  \eqref{charerrorequal} will be solved. The choice of $p$ and $q$ imposes a minimal number  $n_{min} = q-p+1$. 
\item
We determine a vector $\mathbf{\Psi} ^\ast  = \left( \psi_0, \psi_1, ... , \psi_{n-1} \right)' $ that consists of the $n$ first components of the set $\mathbf{\Psi} = \left\{ \psi_0, \psi_1, ... \right\} $ that satisfies condition \eqref{charerrorequal}. We emphasize that in general this condition does not determine uniquely the vector $\mathbf{\Psi} ^\ast $ and that arbitrary choices need to be made. The vector $\mathbf{\Psi*}$ is a truncation at order $n-1$ of the MA representation of the ARMA process that we want to construct. 
\item
We conclude the construction of the ARMA(p,q) model that we are after by designing either an AR(p) polynomial $\boldsymbol{\Phi} $ consistent with causality or a MA(q) polynomial $\boldsymbol{\Theta} $ consistent with invertibility. Then:
\begin{itemize}
\item In the first case, the required model is given by 
\begin{equation}
\label{truncation 1}
\boldsymbol{\Phi}(L) X  = \boldsymbol{\Theta} ^\ast (L) \varepsilon, \quad \mbox{with} \quad \boldsymbol{\Theta}^\ast =\mathbf{\Psi} ^\ast \cdot \boldsymbol{\Phi}.
\end{equation}
\item In the second case, the required model is given by 
\begin{equation}
\label{truncation 2}
\boldsymbol{\Phi}^\ast (L) X  = \boldsymbol{\Theta}  (L) \varepsilon, \quad \mbox{with} \quad \boldsymbol{\Phi}^\ast =(\mathbf{\Psi} ^\ast) ^{-1} \cdot \boldsymbol{\Theta}.
\end{equation}
\end{itemize}
In both cases, the MA and AR polynomials that are obtained in this way have to be checked regarding invertibility and causality, respectively. Additionally, the finite truncation of $\boldsymbol{\Psi}$ is likely to give rise to common roots between the AR and MA polynomials in~(\ref{truncation 1}) or in~(\ref{truncation 2}) which may make necessary a slight perturbation of the coefficients in order to be avoided.
\item
We emphasize that the resulting ARMA model satisfies  \eqref{charerrorequal} only approximately and hence the characteristic errors of the two predictors will be not identical but just close to each other for the specific aggregation period $K$ used. For pure MA models no truncation is necessary and hence exact equality can be achieved.
\end{itemize}

\subsection{Stock aggregation examples}
\label{Stock aggregation examples}
In the particular case of stock temporal aggregation,  condition \eqref{charerrorequal} is written as:
\begin{equation} \label{stockCond}
\mathbf{\Psi} \left( L \right) = \left( \sum^{\infty}_{j = 0} \psi_{jK}  L^{jK} \right) \left( \sum^{K-1}_{j = 0}  \psi_{j}  L^{j} \right).
\end{equation} 
We now consider the truncated vector $\mathbf{\Psi^*}$ with $n$ components, that is, $\mathbf{\Psi^*} = \left( \psi_0, \dots, \psi_{n-1}\right)'  $. Then, the truncated version of  \eqref{stockCond} is:
\begin{equation} \label{stockCond_n}
 \sum^{n-1}_{j = 0} \psi_{j}  L^{j}  = \left( \sum^{\lfloor(n-K+1)/K\rfloor}_{j = 0} \psi_{jK}  L^{jK} \right) \left( \sum^{K-1}_{j = 0}  \psi_{j}  L^{j}\right),
\end{equation} 
where the symbol $\lfloor \cdot \rfloor $ denotes the integer part of its argument. We now provide a few examples of  models whose specification is obtained following the approach proposed in the previous subsection and the relation~(\ref{stockCond_n}).
\begin{example}
{\rm MA(10) model.}
\end{example}
Let $p = 0$, $q = 10$, $n = n_{min} = 11$ and let $K=2$.
Equation \eqref{stockCond_n} becomes
$$
\sum^{10}_{j = 0} \psi_{j}  L^{j} = \left( \sum^{5}_{j = 0} \psi_{2j}  L^{2j} \right) \left( \sum^{1}_{j = 0}  \psi_{j}  L^{j} \right),
$$
which imposes the following relations:
$$
\psi_0 = 1,\enspace \psi_1 \psi_{2i} = \psi_{2i+1}, \enspace i = 0,\dots, 5; \enspace \psi_i = 0 \enspace \text{for} \enspace i\ge n.\\
$$
This system of nonlinear equations has many solutions. We choose one of them by setting $\psi_j = 0$, for $j = 1, \dots, 9$, and $\psi_{10} = 0.3$. This way we can trivially determine a MA(10) model which satisfies exactly the relation~(\ref{stockCond}) by taking $\theta_j = 0$ for $j = 1, \dots, 9$ and $\theta_{10} = 0.3$.

\begin{figure}[htp]
\centering
\includegraphics[scale=0.45]{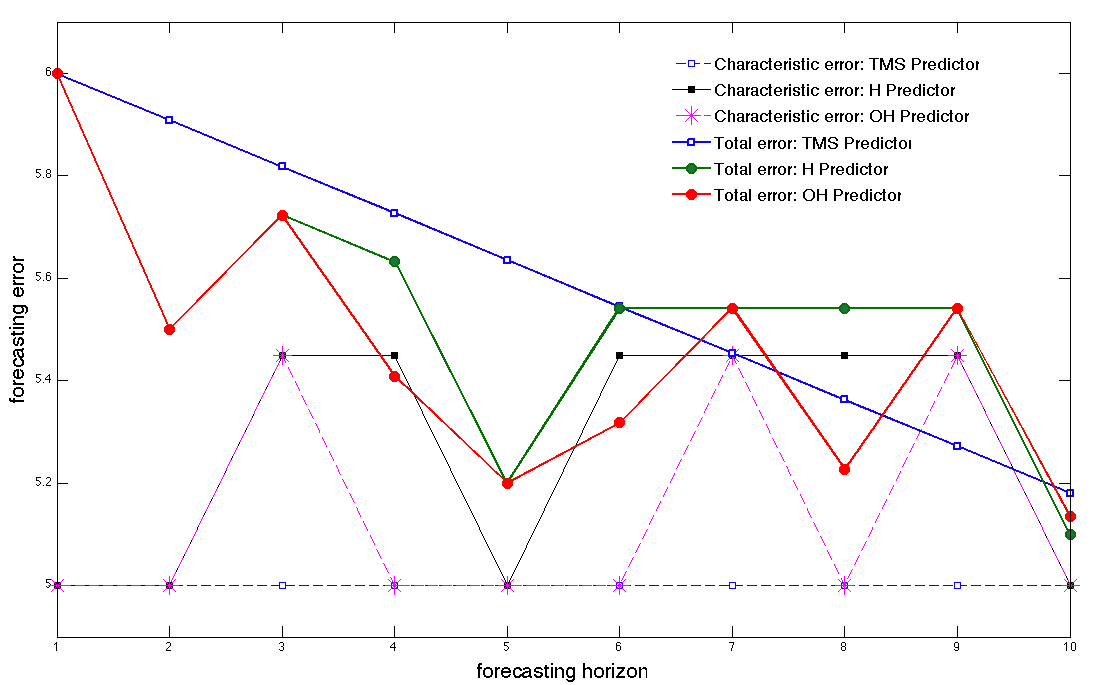}
\caption{Characteristic and total errors associated to the forecast of the temporal stock aggregate of the MA(10) model. The sample size used for estimation is $T = 50$. The innovations of the model have variance $\sigma^2 = 5$.}
\label{fig:MA10_stock}
\end{figure}

Figure~\ref{fig:MA10_stock} shows the values of the characteristic errors for different values of the forecasting horizon for the TMS predictor, the H predictor, and the OH predictor. For the horizon $h = 2$, the values of the characteristic errors of all the predictors coincide, which is a consequence of the fact that the model has been constructed using the relation~(\ref{stockCond}) with $K = 2$. Moreover, it is easy to see by looking at~(\ref{periodic condition}), that the particular choice of MA coefficients that we have adopted ensures that the resulting model is seasonal for the periods 2,  5, and 10; this guarantees that~(\ref{stockCond}) is also satisfied for the corresponding values of $K$ and hence there is coincidence for the characteristic errors at those horizons too.

The total errors for a sample size of $T=50$ are then computed using the formulas presented in sections~\ref{Forecasting with estimated linear processes} and~\ref{Finite sample forecasting of temporally aggregated linear processes}. The corresponding results  are also plotted in Figure~\ref{fig:MA10_stock} for the different forecasting schemes. This plot shows that for several forecasting horizons both the H and the OH predictors perform better than the TMS predictor. 

A quick inspection of this plot reveals another interesting phenomenon consisting on the decrease of the total error associated to the three predictors as the forecasting horizon increases; this feature is due to the decrease of the estimation error using these forecasting schemes as the horizon becomes longer.  The characteristic errors for the H and OH  predictors do not increase monotonically with the forecasting horizon either; however, in this case, this is due to the fact that for each value of the forecasting horizon, these predictors are constructed using a different model since the aggregation period changes and hence so does the aggregated model used for forecasting.

In conclusion, in this particular example, both the H  and the OH predictors exhibit a better forecasting performance than the TMS predictor and, additionally, the results regarding the OH predictor help in making a decision on what is the best possible aggregation period to work with in order to minimize the associated total forecasting error.   

\begin{figure}[htp]
\centering
\includegraphics[scale=0.45]{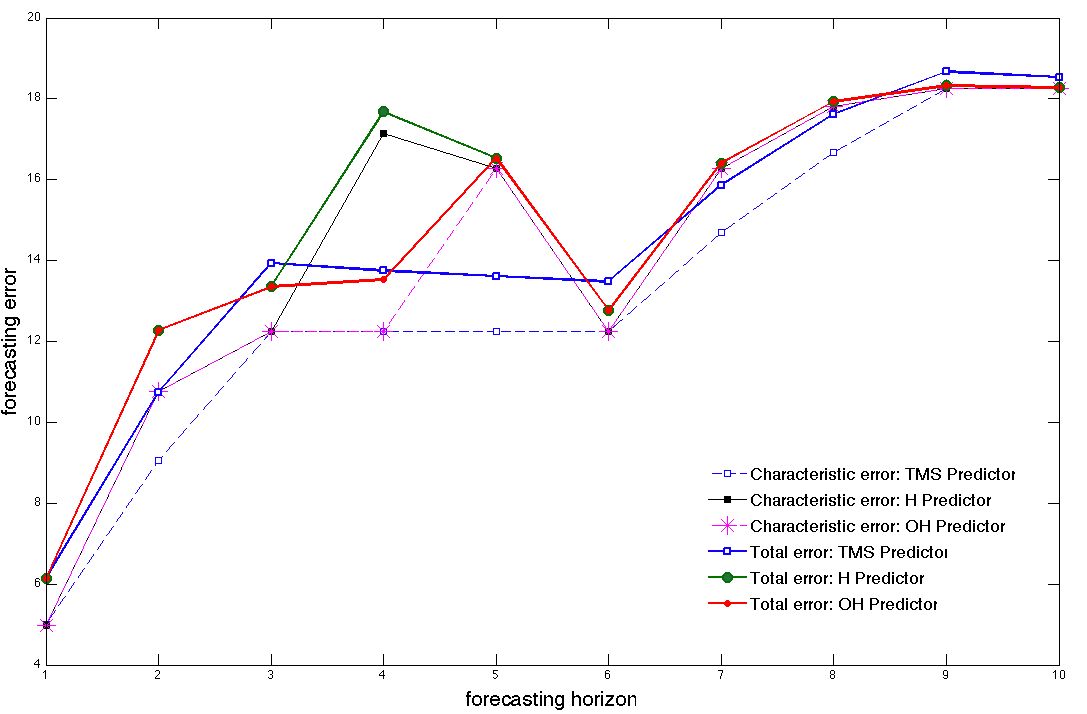}
\caption{Characteristic and total errors associated to the forecast of the temporal stock aggregate of the ARMA(3,11) model. The sample size used for estimation is $T = 50$. The innovations of the model have variance $\sigma^2 = 5$.}
\label{fig:ARMA311_stock}
\end{figure}

\begin{example}
{\rm ARMA(3,11)} model.
\end{example}
Let $p = 3$, $q = 11$, $n = 10$ and let $K=3$.
In this case, the relation~\eqref{stockCond_n} yields
$$
\sum^{9}_{j = 0} \psi_{j}  L^{j} = \left( \sum^{2}_{j = 0} \psi_{3j}  L^{3j} \right) \left( \sum^{2}_{j = 0}  \psi_{j}  L^{j} \right)
$$
and consequently
$$
\psi_0 = 1,\enspace \psi_j \psi_{3i} = \psi_{3i+1}, \enspace i = 0,\dots, 2, \enspace j = 1, 2, \enspace \psi_i = 0 \enspace \text{for} \enspace i \ge n-1.\\
$$
We choose a solution for these relations of the form $\mathbf{\Psi^*} = (1, -0.9, 0.8,$ $0, 0, 0, -0.7,$ $0.63, -0.56, 0)'$. We now introduce an AR(3) polynomial of the form  $\mathbf{\Phi} = (-0.9, 0.8, -0.4)'$. We then determine the MA(11) part of the model by using~(\ref{truncation 1}), which yields the coefficients $\mathbf{\Theta} = (-1.8, 2.41,-1.84,$ $1, -0.32, -0.7, 1.26,$ $-1.687, 1.288, -0.7, 0.224)'$. In order to avoid the common roots between the AR and the MA polynomials that are  obtained when the coefficients of the MA part are derived in this manner, we slightly perturb the values of some of the components of the vector $\mathbf{\Theta}$ that we now set to be $\mathbf{\Theta} = (-1.8, 2.4102,-1.8403, 1, -0.32,$ $ -0.7, 1.26, -1.687,$ $ 1.288, -0.7, 0.224)'$. Figure~\ref{fig:ARMA311_stock} shows the errors with respect to the forecasting horizon for all the predictors as in the previous example. The H and the OH predictors perform better than the TMS predictor for $h=3, 6, 9, 10$. Additionally the OH predictor performs better than the H predictor for $h = 4$. Taking into account the initial choice of $K=3$ when constructing the example, it becomes clear why the characteristic errors associated to  the H and the OH predictors are very close to those associated to the TMS predictor for horizons $h$ that are multiples of $3$.  
\begin{example}
{\rm ARMA(1,4)} model.
\end{example}
Let $p = 1$, $q = 4$, $n = 5$ and let $K=4$.
In this setup, relation \eqref{stockCond_n} yields
$$
\sum^{4}_{j = 0} \psi_{j}  L^{j} = \sum^{3}_{j = 0}  \psi_{j}  L^{j}, 
$$
and consequently
$\psi_0 = 1$ and  $ \psi_4  = 0 $, necessarily, while the values of the coefficients $ \psi_1$, $\psi_2$, and  $\psi_3$ are not subjected to any constraint. We hence set
 $\mathbf{\Psi^*} = (1, 0.3, -0.3, 0.3, 0)'$. 
 We now introduce the AR(1) polynomial determined by the coefficient $ \boldsymbol{\Phi}=0.8 $. We then determine the MA(4) part of the model by using~(\ref{truncation 1}) which yields $\mathbf{\Theta} =  (-0.5, -0.54, 0.54, -0.24)'$. Again in order to avoid  common roots between the AR and the MA polynomials, we perturb the polynomial $\mathbf{\Theta}$ by setting: $\mathbf{\Theta} =  (-0.5, -0.5403, 0.54, -0.24)'$. 
 
\begin{figure}[htp]
\centering
\includegraphics[scale=0.35]{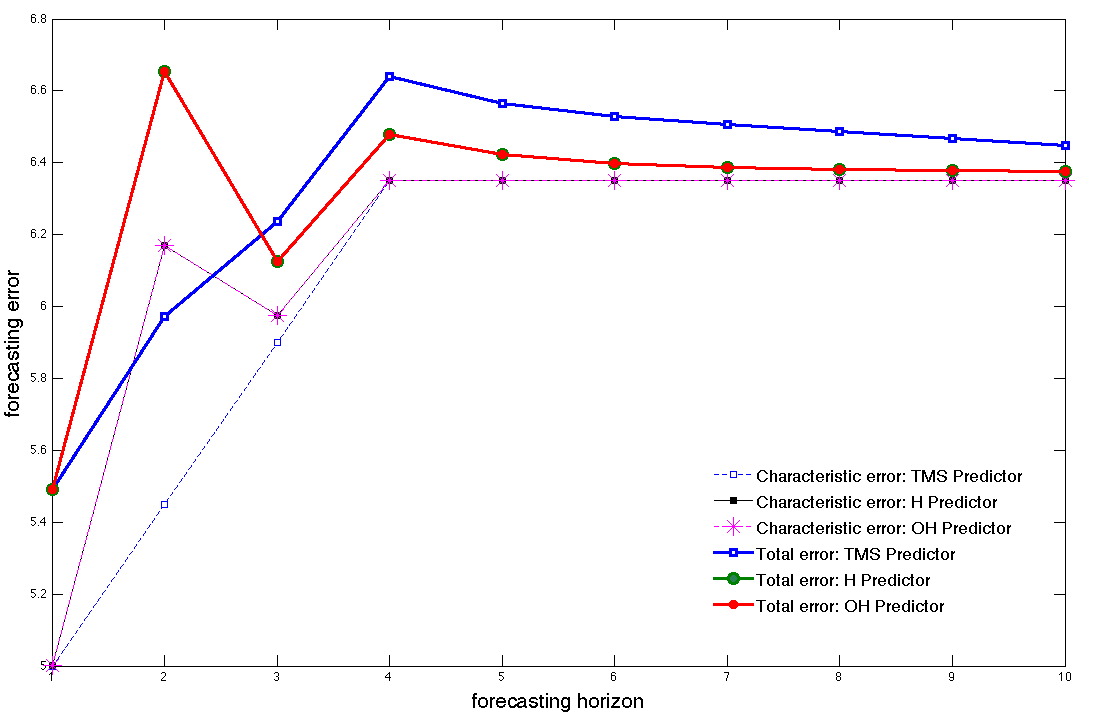}
\caption{Characteristic and total errors associated to the forecast of the temporal stock aggregate of the ARMA(1,4) model. The sample size used for estimation is $T = 50$. The innovations of the model have variance $\sigma^2 = 5$.}
\label{fig:ARMA14_stock}
\end{figure}

Figure~\ref{fig:ARMA14_stock} shows that the H  and the OH predictors have equal associated total errors and exhibit a better forecasting efficiency than the TMS predictor for all  forecasting horizons except at $h=2$. The initial choice of $K=4$ at the model construction stage results in the fact that for $h=4$ the values of the characteristic errors associated to the H and the OH predictors  are very close to the one committed by the TMS predictor.
\subsection{Flow aggregation examples}
\begin{figure}[htp]
\centering
\includegraphics[scale=0.4]{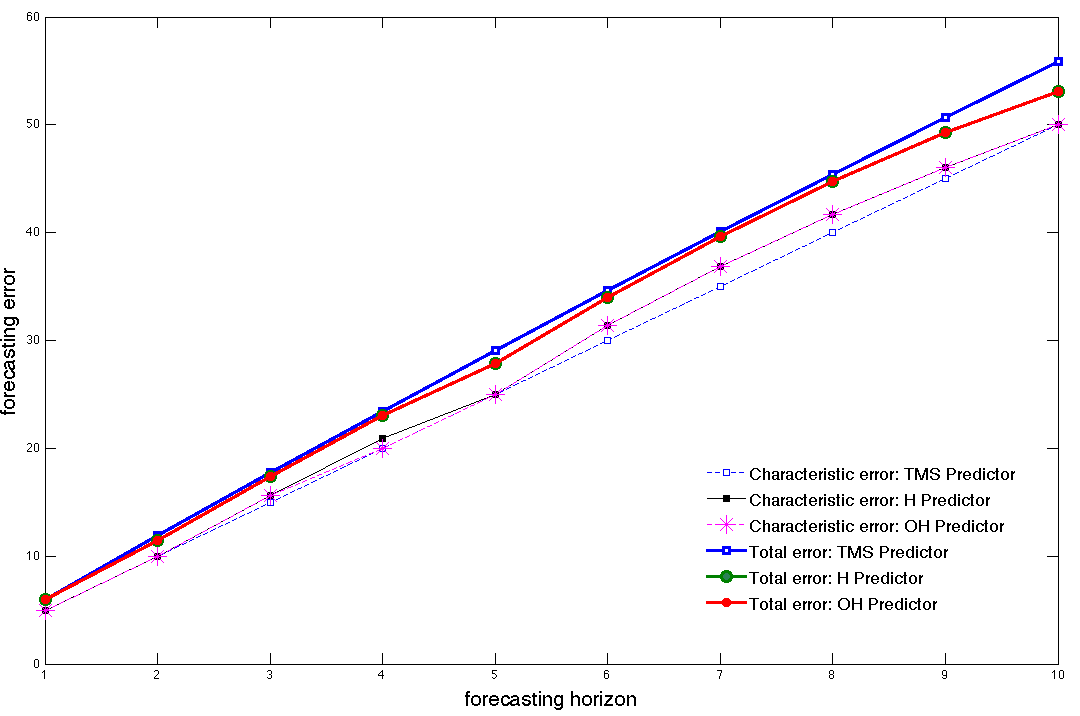}
\caption{Characteristic and total errors associated to the forecast of the temporal flow aggregate of the MA(10) model. The sample size used for estimation is $T = 50$. The innovations of the model have variance $\sigma^2 = 5$.
}
\label{fig:MA10_flow}
\end{figure}

In the particular case of flow temporal aggregation, condition \eqref{charerrorequal} can be written as:
\begin{equation} \label{flowCond}
 \mathbf{\Psi} \left( L \right)\sum^{K-1}_{i = 0} L^i  = \left( \sum^{\infty}_{j = 0} L^{jK}  \sum^{K-1}_{i = 0}   \psi_{jK-i}\right)  \left( \sum^{K-1}_{j = 0} L^{j}  \sum^{j}_{i = 0}  \psi_{j-i}  \right).
\end{equation} 
We now consider the truncation $\mathbf{\Psi}^\ast $of  $\mathbf{\Psi}$ with $n $ components, that is, $\mathbf{\Psi^*} = \left( \psi_0, \dots, \psi_{n-1}\right) $. Then, the truncated version of \eqref{flowCond} can be expressed as:
\begin{equation} \label{flowCond_n}
 \left( \sum^{K-1}_{i = 0} L^i \right)  \sum^{n-1}_{j = 0} \psi_{j}  L^{j} = \left( \sum^{\lfloor (n-K+1)/K\rfloor }_{j = 0} L^{jK}  \sum^{K-1}_{i = 0}  \psi_{jK-i} \right)   \left( \sum^{K-1}_{j = 0} L^{j}  \sum^{j}_{i = 0}  \psi_{j-i}  \right),
\end{equation} 
where symbol $\lfloor \cdot \rfloor $ denotes the integer part of its argument.
We now provide a few examples of  models whose specification is obtained following the approach described in the beginning of the section and the relation~(\ref{flowCond_n}). 
\begin{example}
{\rm MA(10)} model.
\end{example}
Let $p = 0$, $q = 10$, $n = n_{min}=11$, and let $K=2$. Then the expression \eqref{flowCond_n} reads
$$
\left( \sum^{1}_{i = 0} L^i \right)  \sum^{10}_{j = 0} \psi_{j}  L^{j} = \left( \sum^{5}_{j = 0} L^{2j}  \sum^{1}_{i = 0}  \psi_{2j-i} \right)   \left( \sum^{1}_{j = 0} L^{j}  \sum^{j}_{i = 0}  \psi_{j-i}  \right),
$$
and consequently
\begin{equation}
\label{system ma10}
\psi_0 = 1,\enspace (1 + \psi_1)( \psi_{i} + \psi_{i+1})=\psi_{i+1} + \psi_{i+2}, \quad \mbox{for} \quad i = 1,\dots, 9, \quad \mbox{and} \quad \enspace \psi_i = 0 \enspace \text{for} \enspace i\ge n.
\end{equation}
A solution for these equations is given by the choice $\psi_j = 0$, for $j = 1, \dots, 9$, and $\psi_{10} = 0.3$. Since the order of the AR polynomial is zero, the method that we proposed determines uniquely in this case the MA(10) polynomial that we are after with $\mathbf{\Theta} = (0,0,0,0,0,0,0,0,0,0.3)'$. The evolution of the forecasting errors versus the forecasting horizon in plotted in the Figure~\ref{fig:MA10_flow}. Both the H and the OH predictors perform better than the TMS predictor. For $h=4$ the OH predictor has the smallest total error among the three predictors. 

\begin{figure}[htp]
\centering
\includegraphics[scale=0.45]{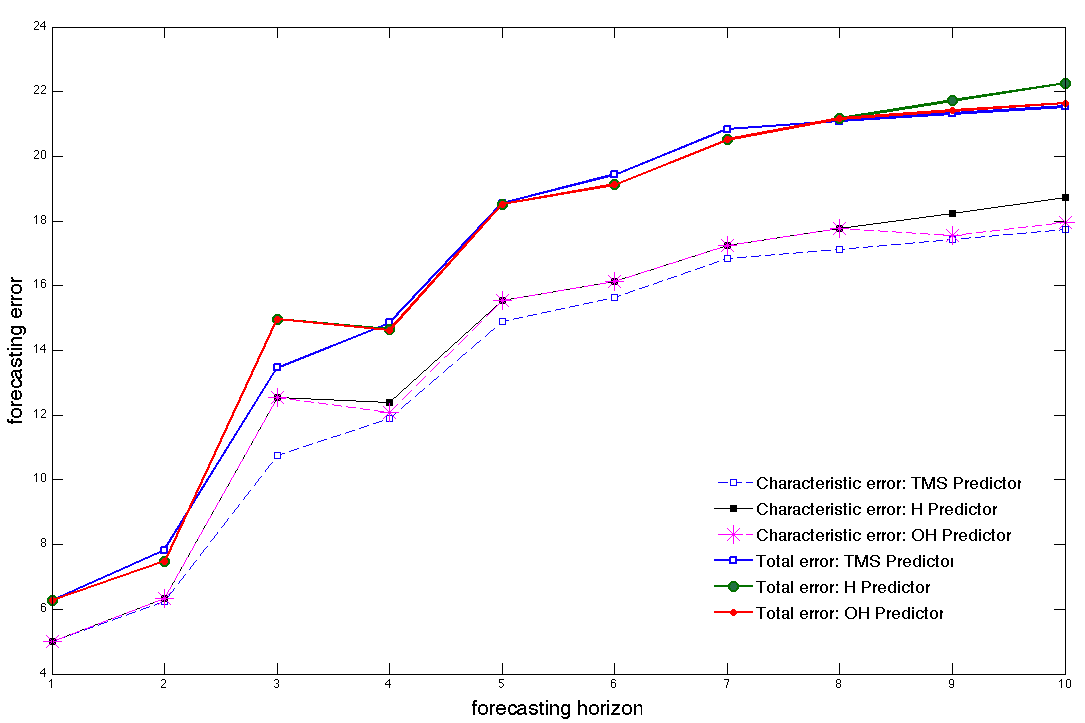}
\caption{Characteristic and total errors associated to the forecast of the temporal flow aggregate of the ARMA(3,10) model. The sample size used for estimation is $T = 50$. The innovations of the model have variance $\sigma^2 = 5$.}
\label{fig:ARMA310_flow}
\end{figure}

\begin{example}
{\rm ARMA(3,10)} model.
\end{example}
In the previous example we chose $\psi_{10} \ne 0$. Let us now use  another solution of the system~(\ref{system ma10}) in order to obtain another model with target orders $p=3$ and $q=10$.
If we set $\psi_{10}=0$, then a possible solution is $\mathbf{\Psi^*} = (1, -0.5, 0.45, -0.475, 0.3, -0.3875, 0.1, -0.2438, 0, 0, 0)'$.  Let now $\mathbf{\Phi} = ( 0.21, 0.207, 0.0162)'$ be a causal AR(3) polynomial which determines via~(\ref{truncation 1}) the MA(10) polynomial  $\mathbf{\Theta} = (-0.71, 0.348, -0.4822, 0.3147,$ $ -0.3595, 0.1270, -0.1894,$ $0.0368, 0.0488, 0.0039)'$. In order to avoid common roots for the AR and MA polynomials, we perturb the MA coefficients and set $\mathbf{\Theta} = (-0.71, 0.3481, -0.4823, 0.3148,$ $ -0.3595, 0.1270,$ $ -0.1894, 0.0368, 0.0488, 0.0039)'$. Figure~\ref{fig:ARMA310_flow} shows the corresponding characteristic and total errors for the three predictors. The H  and the OH predictors exhibit better performance than the TMS predictor for horizons $h = 2,4,5,6, 7$.

\begin{remark}
\normalfont
The model that we just presented can be used to illustrate the fact that the construction method that we presented in this sections is not the unique source of examples in which the H and the OH predictors perform better than the TMS scheme.  Indeed, as it can be seen in Figure~\ref{fig:ARMA310_stock}, the very same ARMA(3,10) model prescription used in the context of stock aggregation also shows this feature even though it has not been obtained by finding a solution of the equation~(\ref{stockCond}).
\begin{figure}[htp]
\centering
\includegraphics[scale=0.45]{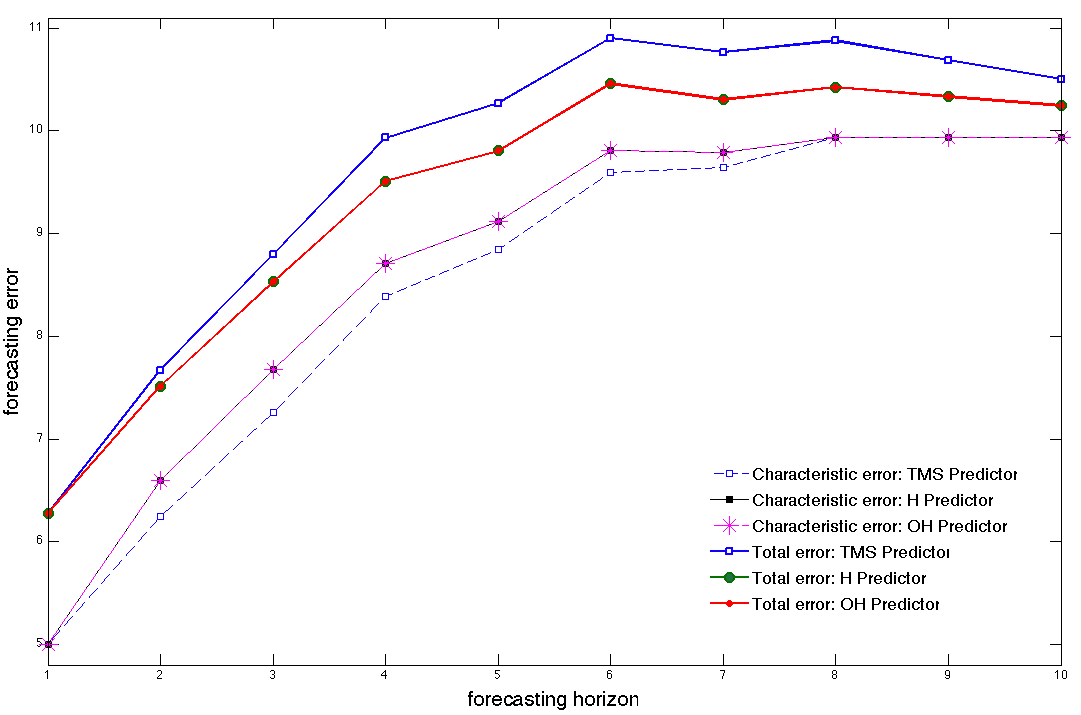}
\caption{Characteristic and total errors associated to the forecast of the temporal stock aggregate of the ARMA(3,10) model. The sample size used for estimation is $T = 50$. The innovations of the model have variance $\sigma^2 = 5$.}
\label{fig:ARMA310_stock}
\end{figure}
\end{remark}

\begin{remark}
\normalfont
Notice that when the forecasting horizon $h$ equals one all predictors coincide because there is no temporal aggregation and hence they obviously have the same errors associated.
\end{remark}

\section{Conclusions}

In this work we have carried out a detailed study of the total error committed when forecasting with one dimensional linear models by minimizing the mean square error. We have introduced a new hybrid scheme for the forecasting of linear temporal aggregates that in some situations shows optimal performance in comparison with other prediction strategies proposed in the literature. 

We work in a finite sample context. More specifically, the forecasting is based on the information set generated by a sample and a model whose parameters have been estimated on it and we avoid the use of second order stationarity hypotheses or the use of time independent autocovariance functions.

In this setup, we provide explicit expressions for the forecasting error that incorporate both the error incurred in due to the stochastic nature of the model (we call it characteristic error) as well as the one associated to the sample based estimation of the model parameters (estimation error). In order to derive these expressions we use certain independence and asymptotic normality hypotheses that are customary in the literature; our main contribution consists of providing expressions for the total error that do not require neither stationarity on the samples used nor Monte Carlo simulations to be evaluated.

We subsequently use these formulas to evaluate the performance of a new forecasting strategy that we propose for the prediction of linear temporal aggregates; we call it {\it hybrid scheme}. This approach consists of using high frequency data for estimation purposes and the corresponding temporally aggregated data and model for forecasting.  This scheme uses all the information available at the time of estimation by using the bigger sample size provided by the disaggregated data, and allows these parameters to be updated as new high frequency data become available. More importantly, as we illustrate with various examples, in some situations the total error committed using this scheme is {\it smaller} than the one associated to the forecast based on the disaggregated data and model; in those cases our strategy becomes optimal. As the increase in performance obtained with our method comes from minimizing the estimation error, we are persuaded that this approach to forecasting may prove very relevant in the multivariate setup where in many cases the estimation error is  the main source of error.

\section{Appendices}

\subsection{Proof of Proposition~\ref{prop1}}

\begin{enumerate}
\item[\bf{(i)}]
It is a straightforward consequence of the causality and invertibility hypotheses on the ARMA model that we are dealing with. Indeed, we can write
\begin{equation} \label{eq:prop11proof}
\varepsilon_t = \sum^{t - 1 + r}_{j = 0}  \pi_j X_{t - j} \quad \text{and} \quad X_t = \sum^{t - 1 + r}_{i = 0} \psi_i  \varepsilon_{t - i},
\end{equation} 
which proves \eqref{eq:prop11}.
\item[\bf{(ii)}]
Suppose first that the innovations $\left\{ \varepsilon_t \right\} $ are IID$(0, \sigma ^2)$. Then the forecast $\widehat{X_{T + h}}$ that minimizes the mean square forecasting error $ E \left[ \left( X _{T+h} - \widehat{ X _{T+h}} \right) ^2 \right]  $ is given by the conditional expectation (see for example~\cite{MR1278033}, page 72):
\begin{align}
\widehat{X_{T + h}} & = E \left[ X_{T + h} | \sigma \left( \underline{\xi_{T}}\right) \right] = E  \left[ X_{T + h} | \sigma \left( \underline{\epsilon_{T}}\right) \right]  = \sum^{T + h - 1 + r}_{i = 0} \psi _{i} E  \left[ \varepsilon _{T + h - i} | \sigma \left( \underline{\epsilon_{T}}\right) \right]  \nonumber \\ 
& = \sum^{T + h -1 + r}_{i = h} \psi_{i} \varepsilon_{T + h - i} = \sum^{T - 1 + r}_{i = 0} \psi_{i + h} \varepsilon_{T - i} = \sum^{T - 1 + r}_{i = 0} \sum^{T - i - 1 + r}_{j = 0} \psi _{i + h} \pi _{j} X_{T - i - j}, \nonumber 
\end{align}
as required.

When $\left\{ \varepsilon_t \right\} $ is WN$(0, \sigma ^2)$ our goal is finding the linear combination $ \sum^{T - 1 + r}_{j = 0} a_{j} X_{T - j} $ that minimizes 
\begin{equation*} 
E \left[ \left( X _{T+h} - \sum^{T - 1 + r}_{i = 0} a_{i} X_{T - i} \right) ^2 \right].
\end{equation*} 
Given that by \eqref{eq:prop11proof}, the elements $ X_{T - i} $ can be written as a linear combination of the elements in $\epsilon_{T - i} $, this task is equivalent to finding the vector $ \mathbf{b} = \left( b_{0}, ..., b_{T - 1 + r} \right)' $ that minimizes the function
\begin{align}
& S \left( b_{0}, ..., b_{T - 1 + r} \right)  = E  \left[ \left( X _{T+h} - \sum^{T - 1 + r}_{i = 0} b_{i} \varepsilon _{T - i} \right) ^2 \right] = E  \left[ \left( \sum^{T + h - 1 + r}_{i = 0} \psi _i \varepsilon _{T+h-i} - \sum^{T - 1 + r}_{i = 0} b_{i} \varepsilon _{T - i} \right) ^2 \right]  \nonumber \\
& = E  \left[ \left( \sum^{h - 1}_{i = 0} \psi _i \varepsilon _{T+h-i} + \sum^{T - 1 + r}_{i = 0} \left( \psi _{i+h} - b_{i} \right) \varepsilon _{T - i} \right) ^2 \right] = \sigma ^2 \left[  \sum^{h - 1}_{i = 0} \psi _i ^2 + \sum^{T - 1 + r}_{i = 0} \left( \psi _{i+h} - b_{i} \right) ^2 \right].  \nonumber 
\end{align}
Hence, in order to minimize the function $S \left( b_{0}, ..., b_{T - 1 + r} \right) $ we compute the partial derivatives $\partial S/\partial  b _i $  and we set them to zero, which shows that the optimal values are attained when $b _i = \psi _{i+h}$.
Consequently, the optimal linear forecast is given by $ \widehat{X_{T + h}} = \sum^{T - 1 + r}_{i = 0} \psi _{i+h} \varepsilon _{T - i} $, as required.
\item[\bf{(iii)}]
We first compute $ X_{T + h} - \widehat{X_{T + h}}$. By \eqref{eq:prop12} and \eqref{eq:prop11proof} we have
$$X_{T + h} - \widehat{X_{T + h}} = \sum^{h-1}_{i= 0}\psi _i \varepsilon _{T+h-i}.$$
Therefore
$$\text{MSFE} \left( \widehat{X_{T + h}} \right) = E  \left[ \left( \sum^{h-1}_{i = 0} \psi _i \varepsilon _{ T+h - i} \right) ^2 \right] = \sigma ^2 \sum^{h-1}_{i = 0} \psi _i ^2.$$
\item[\bf{(iv)}]
Given the model $ \mathbf{\Phi} \left( L \right) X = \mathbf{\Theta} \left( L \right) \varepsilon$, we have 
\begin{equation*} 
X _{ T+h} -\phi_{1} X_{T + h - 1} + ... + \phi_{p} X_{T + h - p} = \varepsilon _{ T+h} +  \theta_{1} \varepsilon_{T+h-1} + ... + \theta_{q} \varepsilon_{T + h - q}.
\end{equation*} 
We first recall that by \eqref{eq:prop11} we have that $\sigma \left( \underline{\xi_T} \right) = \sigma \left( \underline{\epsilon  _T } \right)=:\mathcal{F}_T  $. We now project both sides of this equality onto the information set $\mathcal{F}_T $ by thinking of this $\sigma$-algebra as $\sigma \left( \underline{\xi_T} \right)  $ for the left hand side projection and as $\sigma \left( \underline{\epsilon  _T } \right)  $  for the right hand side.  We obtain:
\begin{align}
\widehat{X_{T + h}} - \phi_{1} \widehat{ X_{T + h - 1}} + ... + \phi_{p} \widehat{ X_{T + h - p}} & = E  \left[ \varepsilon _{ T+h} +  \theta_{1} \varepsilon_{T+h-1} + ... + \theta_{q} \varepsilon_{T + h - q} | \mathcal{F}_T) \right] \nonumber \\
& = \left\{
\begin{array}{l l}
\theta_{h} \varepsilon_{T} + ... + \theta_{q} \varepsilon_{T + h - q}, & q \ge h\\
0 , & \text{otherwise}. \\
\end{array} \right. \nonumber 
\end{align}
In the presence of white noise innovations, the conditional expectation in the previous equality should be replaced by a linear projection. \quad $\blacksquare$
\end{enumerate} 

\subsection{Computation of the Jacobian $J_{\boldsymbol{\Xi}_P}$}
\label{computation of the jacobian}
In this section we provide a simple algorithmic construction for the computation of the Jacobian $J_{\mathbf{\Xi}_P} $ when the elements in $\mathbf{\Xi}_P $ are considered as a function of $\boldsymbol{\beta} $, that is, $\mathbf{\Xi}_P(\boldsymbol{\beta}) := \left( \psi _1(\boldsymbol{\beta}) , \dots, \psi_P(\boldsymbol{\beta}) , \pi _1 (\boldsymbol{\beta}), \dots, \pi _P(\boldsymbol{\beta}) \right)$. We will separately compute the components $\dfrac{\partial{\psi_i}}{\partial{\beta_k}}$ and $\dfrac{\partial{\pi_i}}{\partial{\beta_k}}$, $i=1, \ldots, P $,  $k=1,\dots, p+q$.

The causality and invertibility hypotheses on the ARMA process we are working with, guarantee that for any $P \ge \max \left\{ p,q\right\} $ there exist polynomials $\mathbf{\Psi}_P (z)$, $\boldsymbol{\Pi}_P(z)$  of order $P$ uniquely determined by the relations:
\begin{align}
\label{psieq}
&\mathbf{\Phi} (z) \mathbf{\Psi}_P (z) = \mathbf{\Theta} (z), \\
\label{pieq}
&\mathbf{\Phi} (z) =  \boldsymbol{\Pi}_P (z) \mathbf{\Theta} (z),
\end{align}
which are equivalent to
\begin{align}
&\mathbf{\Psi}_P (z) = \mathbf{\Phi} ^{-1} (z) \mathbf{\Theta} (z), \nonumber \\
&\boldsymbol{\Pi}_P (z) = \mathbf{\Phi} (z) \mathbf{\Theta} ^{-1} (z). \nonumber
\end{align}
These polynomial relations determine the functions $\mathbf{\Psi}_P ( \mathbf{\Phi}, \mathbf{\Theta})$, $\boldsymbol{\Pi}_P ( \mathbf{\Phi}, \mathbf{\Theta})$ needed in the computation of the Jacobian. We now rewrite \eqref{psieq} and \eqref{pieq} as
\begin{align}
\label{psieqfunc}
& \mathbf{\Phi} (z) \mathbf{\Psi}_P ( \mathbf{\Phi}, \mathbf{\Theta}) (z) = \mathbf{\Theta} (z),  \\
\label{pieqfunc}
& \mathbf{\Phi} (z) = \boldsymbol{\Pi}_P ( \mathbf{\Phi}, \mathbf{\Theta}) (z) \mathbf{\Theta} (z).
\end{align}
If we take derivatives with respect to $\theta _j $ and $\phi _i $, $j \in \left\{ 1, \dots, q\right\} $, $i \in \left\{ 1, \dots, p\right\} $ on both sides of \eqref{psieqfunc},  we obtain a set of $p+q$ polynomial equations:
\begin{align}
& \mathbf{\Phi} (z) \dfrac{\partial{\mathbf{\Psi}_P (\mathbf{\phi}, \mathbf{\theta })}}{\partial{\theta _j }}(z) = z _j, \enspace j \in \left\{ 1, \dots, q \right\}, \nonumber \\
& z ^i \mathbf{\Psi} _P (\mathbf{\phi}, \mathbf{\theta})+\mathbf{\Phi} (z) \dfrac{\partial{\mathbf{\Psi} _P (\mathbf{\phi}, \mathbf{\theta })}}{\partial{\phi _i }} = 0, \enspace i \in \left\{ 1, \dots, p \right\},\nonumber
\end{align}
that determine uniquely the corresponding entries of the Jacobian due to the invertibility of $\mathbf{\Phi} (z)$. At the same time, taking derivatives on both the right and left hand sides of \eqref{pieqfunc} with respect to $\theta _j $ and $\phi _i $, we obtain another set of $p+q$ polynomial equations
\begin{align}
& \dfrac{\partial{\mathbf{\Phi}}}{\partial{\theta_j}} (z) = \dfrac{\partial{\boldsymbol{\Pi}_P (\mathbf{\phi}, \mathbf{\theta })}}{\partial{\theta _j }}(z) + \boldsymbol{\Pi}_P  z _j, \enspace j \in \left\{ 1, \dots, q \right\}, \nonumber \\
& z ^i =   \dfrac{\partial{\boldsymbol{\Pi}_P (\mathbf{\phi}, \mathbf{\theta })}}{\partial{\phi _i }}, \enspace i \in \left\{ 1, \dots, p \right\},\nonumber
\end{align}
that determine uniquely the corresponding entries of the Jacobian due to the invertibility of $\mathbf{\Theta} (z)$.

\subsection{Proof of Theorem~\ref{theor1} } 
\begin{enumerate}
\item[\bf{(i)}]
It is a straightforward consequence of the independence hypothesis between the samples $\xi_T$ and the $\xi'_T$, and part {\bf (ii)} of Proposition~\ref{prop1}.
\item[\bf{(ii)}]
By \eqref{31} and part {\bf (i)} of Proposition~\ref{prop1} we have that
\begin{equation*}
 \text{MSFE} \left( \widehat{X_{T + h}} \right) = E \left[ \left( X _{t+h} - \widehat{ X _{t+h}} \right) ^2 \right]  = E \left[ \left(\sum^{P}_{i = 0} \psi_{i} \varepsilon_{T + h - i} - \sum^{P}_{i = h} \widehat{\psi}_{i} \tilde{\varepsilon}_{T + h - i}\right) ^2 \right].
 \end{equation*}
In order to compute this error notice that  $\tilde{\varepsilon }_{T+h-i} $ can be rewritten in terms of the original innovations as
\begin{equation} \label{prop3:proof:2} 
\tilde{\varepsilon }_{T+h-i} = \sum^{P-i}_{j = 0} \widehat{\pi}_j x_{T+h-i-j} = \sum^{P-i}_{j = 0} \sum^{P-i-j}_{k = 0} \widehat{\pi}_j \psi _k \varepsilon _{T+h-i-j-k}.
\end{equation} 
Hence,
\begin{align} 
\label{prop3:proof:2} 
 & \text{MSFE} \left( \widehat{X_{T + h}} \right) = E \left[ \left(\sum^{P}_{i = 0} \psi_{i} \varepsilon_{T + h - i} - \sum^{P}_{i = h} \sum^{P-i}_{j = 0} \sum^{P-i-j}_{k = 0}  \widehat{\psi}_{i} \widehat{\pi}_j \psi _k \varepsilon_{T + h - i-j-k}\right) ^2 \right]  \nonumber \\  
 & = E  \Biggr[\left(\sum^{P}_{i = 0} \psi_{i} \varepsilon_{T + h - i} \right) ^2 - 2 \sum^{P}_{l = 0} \sum^{P}_{i = h} \sum^{P-i}_{j = 0} \sum^{P-i-j}_{k = 0} \psi _l  \widehat{\psi}_{i} \widehat{\pi}_j \psi _k \varepsilon _{T+h-l} \varepsilon_{T + h - i-j-k} \nonumber \\
 &+ \left.\left( \sum^{P}_{i = h} \sum^{P-i}_{j = 0} \sum^{P-i-j}_{k = 0}  \widehat{\psi}_{i} \widehat{\pi}_j \psi _k \varepsilon_{T + h - i-j-k}\right)^2  \right]  = \sigma ^2 \Biggl[ \sum^{P}_{i = h} {\psi_{i}} ^{2}  - 2 \sum^{P}_{l = 0} \sum^{P}_{i = h} \sum^{P-i}_{j = 0} \sum^{P-i-j}_{k = 0}  \psi_{l} \psi_k E \left[ \widehat{ \psi _i } \widehat{\pi }_j \right]  \delta _{l,i+j+k} \nonumber \\
 &+\sum^{P}_{i = h} \sum^{P-i}_{j = 0} \sum^{P-i-j}_{k = 0} \sum^{P}_{i' = h} \sum^{P-i'}_{j' = 0} \sum^{P-i'-j'}_{k' = 0} \psi _k \psi _{k'} E \left[ \widehat{ \psi} _i \widehat{\pi }_j  \widehat{\psi} _{i'} \widehat{\pi }_{j'}\right]  \delta _{i+j+k,i'+j'+k'}  \Biggr] + \sigma ^2  \sum^{h - 1}_{i = 0} {\psi_{i}} ^{2}. \nonumber
 \end{align} 
\item[\bf{(iii)}]
By part {\bf (i)} the forecast $ \widehat{X_{T + h}}$ is given by
\begin{equation}
\label{recall forecast}
\widehat{X_{T + h}} = \sum^{P}_{i= h}\widehat{\psi} _i \tilde{\varepsilon} _{T+h-i}.
\end{equation}
According to the statement \eqref{eq:theor33}, both $ \widehat{\psi} _i$ and $ \widehat{\pi}_j $ can be asymptotically written as  
$$\widehat{\psi }_i = \psi_i + \dfrac{r_i}{\sqrt{T}}, \enspace \widehat{\pi }_j = \pi_j + \dfrac{t_j}{\sqrt{T}},$$
with $ r_i$ and $t_j$ as Gaussian random variables of mean 0 and variances $(\Sigma _{\mathbf{\Xi}})_{i,i}$ and $(\Sigma _{\mathbf{\Xi}})_{j+P,j+P}$, respectively. Consequently by \eqref{recall forecast} and \eqref{prop3:proof:2}  
\begin{align} 
& \widehat{ X_{T + h}} = \sum^{P}_{i = h} \widehat{ \psi_{i} } {\tilde{\varepsilon}}_{T + h - i} = \sum^{P}_{i = h} \sum^{P-i}_{j = 0} \sum^{P-i-j}_{k = 0} \left(\psi_i + \dfrac{r_i}{\sqrt{T}} \right) \left( \pi_j + \dfrac{t_j}{\sqrt{T}}\right) \psi _k \varepsilon _{T+h-i-j-k} \nonumber \\
& = \sum^{P}_{i = h} \sum^{P-i}_{j = 0} \sum^{P-i-j}_{k = 0} \left(\psi_i  \pi_j \psi _k+ \dfrac{\psi_i t_j  \psi _k}{\sqrt{T}} + \dfrac{r_i  \pi_j \psi _k}{\sqrt{T}} + \dfrac{r_i t_j  \psi _k}{T}\right) \varepsilon _{T+h-i-j-k}.\nonumber
\end{align} 
We now recall that 
\begin{equation*}
\sum^{P-i}_{j = 0} \sum^{P-i-j}_{k = 0} \pi_j \psi _k \varepsilon _{T+h-i-j-k} = \varepsilon _{T+h-i},
\end{equation*}
and we eliminate in this expression the term that decays as 1/T; we hence approximate $ \widehat{X_{T+h}}$ as
\begin{align}
\widehat{X_{T + h}} &\simeq\sum^{P}_{i = h} \left[ \left( \psi_i + \dfrac{r_i}{\sqrt{T}}  \right) \varepsilon _{T+h-i} + \sum^{P-i}_{j = 0} \sum^{P-i-j}_{k = 0} \dfrac{\psi _i \psi _k t _j }{\sqrt{T}} \varepsilon _{T+h-i-j-k}\right]  \nonumber \\
& = \sum^{P}_{i = h}  \widehat{\psi}_i \varepsilon _{T+h-i} + \sum^{P}_{i = h}  \sum^{P-i}_{j = 0} \sum^{P-i-j}_{k = 0} \dfrac{\psi _i \psi _k t _j }{\sqrt{T}} \varepsilon _{T+h-i-j-k}. \nonumber 
\end{align}
Using this approximation we compute now the MSFE:
\begin{align}
\text{MSFE} \left( \widehat{X_{T + h}} \right) &= E \left[ \left( X _{t+h} - \widehat{ X _{t+h}} \right) ^2 \right]   \nonumber \\
 & =  E \left[ \left(\sum^{P}_{i = 0} \psi_{i} \varepsilon_{T + h - i} - \sum^{P}_{i = h}  \widehat{\psi}_i \varepsilon _{T+h-i} - \sum^{P}_{i = h}  \sum^{P-i}_{j = 0} \sum^{P-i-j}_{k = 0} \dfrac{\psi _i \psi _k t _j }{\sqrt{T}} \varepsilon _{T+h-i-j-k}  \right) ^2 \right]  \nonumber \\
 & = E \left[ \left(\sum^{h-1}_{i = 0} \psi_{i} \varepsilon_{T + h - i} + \sum^{P}_{i = h}  \left( \psi _i - \widehat{\psi}_i \right) \varepsilon _{T+h-i} - \sum^{P}_{i = h}  \sum^{P-i}_{j = 0} \sum^{P-i-j}_{k = 0} \dfrac{\psi _i \psi _k t _j }{\sqrt{T}} \varepsilon _{T+h-i-j-k}  \right) ^2 \right] \nonumber \\
 & = \sigma ^2  \sum^{h - 1}_{i = 0} {\psi_{i}} ^{2} 
+ \sigma ^2 \dfrac{1}{T}\Biggl[  \sum^{P}_{i = h}(\Sigma _{\mathbf{\Xi}} )_{i,i}  + 2 \sum^{P}_{i = h} \sum^{P-i}_{j = 0} \sum^{P-i-j}_{k = 0}  \psi_{i} \psi_k  (\Sigma _{\mathbf{\Xi}})_{i+j+k,j+P}  \nonumber \\
& + \sum^{P}_{i = h} \sum^{P-i}_{j = 0} \sum^{P-i-j}_{k = 0} \sum^{P}_{i' = h} \sum^{P-i'}_{j' = 0} \sum^{P-i'-j'}_{k' = 0}  \psi _i \psi _k \psi _{i'} \psi _{k'}  (\Sigma _{\mathbf{\Xi}})_{j+P,j'+P} \delta _{i+j+k,i'+j'+k'} \Biggr]. \quad \blacksquare\nonumber
\end{align}
\end{enumerate}
\subsection{Proof of Theorem~\ref{luetkepohl formula theorem}}

The mean square error associated to the forecast $\widehat{X_{T+h} }( \widehat{ \boldsymbol{\beta} }) $ carried out using  estimated parameters $\widehat{\boldsymbol{\beta}} $  is given by:
\begin{align} 
\label{msfe}
\text{MSFE} \left( \widehat{X_{T+h} }( \widehat{ \boldsymbol{\beta} }) \right) &= E \Bigl[ \left( X_{T+h} - \widehat{X_{T+h} }( \widehat{ \boldsymbol{\beta} }) \right) ^2\Bigr]= E \Bigl[ \left( X_{T+h} -  \widehat{X_{T+h} }( \boldsymbol{\beta} ) +  \widehat{X_{T+h} }(  \boldsymbol{\beta} ) - \widehat{X_{T+h} }( \widehat{ \boldsymbol{\beta} }) \right) ^2\Bigr]  \nonumber \\
& = E \Bigl[ \left( X_{T+h} -  \widehat{X_{T+h} }( \boldsymbol{\beta} ) \right) ^2 \Bigr] +2 E \Bigl[ \left( X_{T+h} -  \widehat{X_{T+h} }( \boldsymbol{\beta} ) \right) \left( \widehat{X_{T+h} }(  \boldsymbol{\beta} ) - \widehat{X_{T+h} }( \widehat{ \boldsymbol{\beta} }) \right) \Bigr] \nonumber \\
&\ \ \ + E \left[  \left( \widehat{X_{T+h} }(  \boldsymbol{\beta} ) - \widehat{X_{T+h} }( \widehat{ \boldsymbol{\beta} }) \right) ^2\right] .
\end{align}
We now recall that 
\begin{equation*}
\widehat{X_{T+h}} \left( \boldsymbol{\beta} \right) = \sum^{P}_{i = h}  \psi _i \varepsilon _{T+h-i}, \quad \mbox{with} \quad P = T+h-1+r,
\end{equation*}
and notice that
\begin{align}
& E \Bigl[ \left( X_{T+h} -  \widehat{X_{T+h} }( \boldsymbol{\beta} ) \right) \left( \widehat{X_{T+h} }\left(  \boldsymbol{\beta} \right) - \widehat{X_{T+h} }( \widehat{ \boldsymbol{\beta} }) \right) \Bigr] = E \Biggl[ \left( \sum^{h-1}_{j = 0} \psi _j \varepsilon _{T+h-j} \right) \left( \sum^{P}_{i = h} \left( \psi _i \varepsilon _{T+h-i} - \widehat{\psi _i} \tilde{\varepsilon}_{T+h-i} \right) \right)  \Biggr] \nonumber \\
& = E \left[ \left( \sum^{h-1}_{j = 0} \psi _j \varepsilon _{T+h-j} \right) \left( \sum^{P}_{i = h} \left( \psi _i \varepsilon _{T+h-i} - \sum^{P-i}_{j = 0} \sum^{P-i-j}_{k = 0} \widehat{\psi _i} \widehat{\pi _j} \psi _k \varepsilon_{T+h-i-j-k} \right) \right)  \right] = 0,\nonumber
\end{align}
since the first term in the product involves the innovations $ \left\{ \varepsilon _{T+1}, \dots, \varepsilon _{T+h} \right\} $ and the second one $ \left\{ \varepsilon _{1-r}, \dots, \varepsilon _{T} \right\} $; these two sets are disjoint and hence independent. Consequently by~(\ref{eq:prop13}) and \eqref{msfe} we have
\begin{equation*}
\text{MSFE} \left( \widehat{X_{T+h} }( \widehat{ \boldsymbol{\beta} }) \right)  = \sigma^2 \sum^{h-1}_{i = 0} \psi^2 _i  + E \left[ \left( \widehat{X_{T+h} }( \boldsymbol{\beta} )  -  \widehat{X_{T+h} }( \widehat{\boldsymbol{\beta}} ) \right) ^2 \right].
\end{equation*}
The second summand of this expression can be asymptotically evaluated using~(\ref{delta use luetkepohl}). Indeed,
\begin{align}
E \left[ \left( \widehat{X_{T+h} }( \boldsymbol{\beta} )  -  \widehat{X_{T+h} }( \widehat{\boldsymbol{\beta}} ) \right) ^2 \right] &= E\left[E \left[ \left( \widehat{X_{T+h} }( \boldsymbol{\beta} )  -  \widehat{X_{T+h} }( \widehat{\boldsymbol{\beta}} ) \right) ^2 | \mathcal{F}_T \right]\right]\notag \\
	&= \dfrac{1}{T} E \left[  \Omega(h) \right] = \dfrac{1}{T}E \left[ \sum^{p+q}_{i,j = 1}J _i J _j \left( \Sigma _{\boldsymbol{\beta}} \right)_{ij}  \right], \end{align}
with 
\begin{equation*}
J _j	:=\frac{ \partial  \widehat{X_{T+h} }( \boldsymbol{\beta} )}{\partial \beta _j}=\frac{ \partial}{\partial \beta _j}  \left( \sum_{k=h}^P\sum_{l=0}^{P-k} \psi_k \pi_l X_{T+h-k-l} \right) =\sum_{k=h}^P\sum_{l=0}^{P-k}
\left(\frac{\partial{\psi_k}}{\partial \beta _j} \pi_l+\psi_k\frac{\partial \pi_l}{\partial \beta _j}\right)X_{T+h-k-l}.
\end{equation*}
Consequently, 
\begin{align}
E \left[  \Omega(h) \right] &= \sum^{p+q}_{i,j = 1}E \left[J _i J _j \left( \Sigma _{\boldsymbol{\beta}} \right)_{ij}  \right]\notag\\
	&=
\sum^{p+q}_{i,j = 1}\sum_{k=h}^P\sum_{l=0}^{P-k}\sum_{m=h}^P\sum_{n=0}^{P-m}
\left(\frac{\partial{\psi_k}}{\partial \beta _j} \pi_l+\psi_k\frac{\partial \pi_l}{\partial \beta _j}\right)\left(\frac{\partial{\psi_m}}{\partial \beta _i} \pi_n+\psi_m\frac{\partial \pi_n}{\partial \beta _i}\right) \left( \Sigma _{\boldsymbol{\beta}} \right)_{ij}E\left[X_{T+h-k-l}X_{T+h-m-n}\right].\label{intermediate exp1}
\end{align}
In the presence of the stationarity hypothesis in part {\bf (ii)} of the theorem we have that 
\[
E\left[X_{T+h-k-l}X_{T+h-m-n}\right]= \gamma(k+l-m-n)
\]
and hence~(\ref{second expresion for luetkepohl}) follows. Otherwise, since we have in general that
\begin{equation}
\label{nonstatapprtheorem}
E\left[ X_{t}X_{s}\right] = \sigma ^2 \sum^{t+r-1}_{i = 0} \sum^{s+r-1}_{j = 0} \psi _i \psi _j \delta _{t-i,s-j},
\end{equation}
we can insert this expression in~(\ref{intermediate exp1}) and we obtain
\begin{align}
E \left[  \Omega(h) \right] &= \sigma^2
\sum^{p+q}_{i,j = 1}\sum_{k=h}^P\sum_{l=0}^{P-k}\sum_{u=0}^{P-k-l}\sum_{m=h}^P\sum_{n=0}^{P-m}\sum_{v=0}^{P-m-n}
\left(
\frac{\partial{\psi_k}}{\partial \beta _j} \frac{\partial{\psi_m}}{\partial \beta _i} \pi_l\pi_n+
\frac{\partial{\psi_k}}{\partial \beta _j} \frac{\partial{\pi_n}}{\partial \beta _i} \pi_l\psi_m+
\frac{\partial{\pi_l}}{\partial \beta _j} \frac{\partial{\psi_m}}{\partial \beta _i} \psi_k\pi_n\right.\notag\\
	&\ \ \ \left.+
\frac{\partial{\pi_l}}{\partial \beta _j} \frac{\partial{\pi_n}}{\partial \beta _i} \psi_k\psi_m\right)
 \left( \Sigma _{\boldsymbol{\beta}} \right)_{ij}\psi _u\psi_v\delta_{k+l+u, m+n+v}\notag\\
 	&= \sigma^2
\sum_{k=h}^P\sum_{l=0}^{P-k}\sum_{u=0}^{P-k-l}\sum_{m=h}^P\sum_{n=0}^{P-m}\sum_{v=0}^{P-m-n}	\left(
(\Sigma _{\mathbf{\Xi}_P})_{k,m}\pi _l\pi _n+(\Sigma _{\mathbf{\Xi}_P})_{k,n+p}\pi_l\psi_m\right.\notag\\
	&\left.\ \ \ +(\Sigma _{\mathbf{\Xi}_P})_{m,l+p}\psi_k \pi _n+(\Sigma _{\mathbf{\Xi}_P})_{l+p,n+p}\psi_k\psi_m
\right)\psi _u\psi_v\delta_{k+l+u, m+n+v}.\label{intermediate exp2}
\end{align}
The required identity~(\ref{first expresion for luetkepohl}) follows directly from~(\ref{intermediate exp2}) by noticing that
\begin{equation*}
\sum_{l=0}^{P-k}\sum_{u=0}^{P-k-l}\left(\pi_l\psi _{u}\delta_{k+l+u, m+n+v}\right)
=\delta_{k, m+n+v},
\end{equation*}
and
\begin{equation*}
\sum_{n=0}^{P-m}\sum_{v=0}^{P-m-n}\left[
\sum_{l=0}^{P-k}\sum_{u=0}^{P-k-l}\left(\pi_l\psi _{u}\delta_{k+l+u, m+n+v}\right)
\right]\pi _n \psi_v=\sum_{n=0}^{P-m}\sum_{v=0}^{P-m-n}
\pi _n \psi_v \delta_{k, m+n+v}= \delta _{k,m}.
\end{equation*}

\subsection{Proof of Proposition~\ref{prop31}} 
\normalfont
First, we notice that by \eqref{TAseriesComp} we have 
\begin{equation*}
\mathcal{F} _{ M} \subset  \mathcal{F} _{ T}.
\end{equation*}
The same relation guarantees that $ X _{ T+k} ^\mathbf{w} = Y _{ M+1}$. Hence the result is a consequence of the following general fact:
\begin{lemma}
\label{lemma1}
Let $ z $ be a random variable in the probability space $ \left( \Omega, P, \mathcal{F} \right) $. Let $ \mathcal{F}^* $ be a sub-sigma algebra of $ \mathcal{F} $ , that is, $ \mathcal{F} ^* \subset  \mathcal{F}  $. Then 
\begin{equation}
\label{lemma1bis} 
E \left[ \left( z - E \left[ z | \mathcal{F} ^* \right] \right) ^2 \right] \ge E \left[ \left( z - E \left[ z | \mathcal{F} \right] \right) ^2 \right].
\end{equation}
\end{lemma}
\noindent\textbf{Proof of Lemma~\ref{lemma1}\ \ } 
\normalfont
\begin{align}
E \left[ \left( z - E \left[ z | \mathcal{F} ^* \right] \right) ^2 \right] & = E \left[ \left( z - E \left[ z | \mathcal{F} ^* \right] - E \left[ z | \mathcal{F} \right] + E \left[ z | \mathcal{F} \right] \right) ^2 \right] = \nonumber E \left[ \left( z - E \left[ z | \mathcal{F} \right] \right) ^2 \right] \\ &+E \left[ \left( E \left[ z | \mathcal{F} \right] - E \left[ z | \mathcal{F} ^* \right] \right) ^2 \right] + 2 E \left[ \left( z - E \left[ z | \mathcal{F} \right] \right)  \left( E \left[ z | \mathcal{F} \right] - E \left[ z | \mathcal{F} ^* \right] \right)\right]. \nonumber
\end{align}
Given that $ E \left[ \left( z - E \left[ z | \mathcal{F} ^* \right] \right) ^2 \right] \ge 0 $, the inequality \eqref{lemma1bis} follows if we show that 
\begin{equation*}
E \left[ \left( z - E \left[ z | \mathcal{F} \right] \right)  \left( E \left[ z | \mathcal{F} \right] - E \left[ z | \mathcal{F} ^* \right] \right)\right] = 0.
\end{equation*}
Indeed, 
\begin{align*}
& E \left[ \left( z - E \left[ z | \mathcal{F} \right] \right)  \left( E \left[ z | \mathcal{F} \right] - E \left[ z | \mathcal{F} ^* \right] \right) | \mathcal{F}\right] = E \left[  z E \left[ z | \mathcal{F} \right] - z E \left[ z | \mathcal{F} ^*\right] - E \left[ z | \mathcal{F} \right] ^2 +  E \left[ z | \mathcal{F} \right] E \left[ z | \mathcal{F} ^* \right] | \mathcal{F} \right] \nonumber \\ & = E \left[ z | \mathcal{F} \right] ^2 - E \left[ z | \mathcal{F} \right] E \left[ z | \mathcal{F} ^* \right] - E \left[ z | \mathcal{F} \right] ^2 + E \left[ z | \mathcal{F} \right]  E \left[ z | \mathcal{F} ^* \right] = 0.  \quad \blacksquare
\end{align*}

\subsection{Proof of Proposition~\ref{prop_char_er_aggr}} 
Part {\bf (i)} is a straightforward consequence of \eqref{eq:prop12}. Regarding {\bf (ii)}, we first have that
\begin{equation*}
 X _{ T+K} ^\mathbf{w} -\widehat{ X _{ T+K} ^\mathbf{w}} = \sum^{K}_{i = 1} w _i \sum^{i-1}_{j = 0} \psi _j \varepsilon _{ T+i-j}.
\end{equation*}
Consequently,
\begin{align}
{\rm MSFE} \left( \widehat{X _{ T+K} ^\mathbf{w}} \right)  & = E \left[ \sum^{K}_{i = 1}  \sum^{K}_{j = 1} w _i w _j \sum^{i-1}_{l = 0} \sum^{j-1}_{m = 0} \psi _l \psi _m \varepsilon _{ T+i-l} \varepsilon _{ T+j-m} \right] \nonumber \\ & = \sigma ^2 \left[ \sum^{K}_{i = 1} w _i ^2  \sum^{i-1}_{l = 0} \psi _l ^2 + 2 \sum^{K-1}_{i = 1} \sum^{K}_{j = i + 1} w _i w _j  \sum^{i-1}_{l = 0} \psi _l  \psi _{ j-i+l} \right].\quad \blacksquare \nonumber 
\end{align}

\subsection{Proof of Theorem~\ref{theor2}} 
\begin{enumerate}
\item[\bf{(i)}] It is a straightforward consequence of part {\bf (i)} in  Theorem~\ref{theor1}.
\item[\bf{(ii)}]
By \eqref{eq:theor41} we have that
\begin{equation*}
 \text{MSFE} \left( \widehat{X_{T + K}^\mathbf{w}} \right) = E \left[ \left( X _{t+K}^\mathbf{w} - \widehat{ X _{t+K}^\mathbf{w}} \right) ^2 \right]  = E \left[ \left(\sum^{K}_{h = 1} w_h\left[ \sum^{P(h)}_{i = 0} \psi_{i} \varepsilon_{T + h - i} - \sum^{P(h)}_{i = h} \widehat{\psi}_{i} \tilde{\varepsilon}_{T + h - i}\right] \right) ^2 \right],
 \end{equation*}
 where $P(h) = T+h-1+r.$
We now notice that
\begin{equation} \label{theor2:proof:21} 
\tilde{\varepsilon }_{T+h-i} = \sum^{P(h)-i}_{j = 0} \sum^{P(h)-i-j}_{k = 0} \widehat{\pi}_j \psi _k \varepsilon _{T+h-i-j-k}.
\end{equation} 
Hence,
\begin{align} 
\label{theor2:proof:22} 
 & \text{MSFE} \left( \widehat{X_{T + K}^\mathbf{w}} \right) = E \left[ \left(\sum^{K}_{h = 1} w_h\sum^{P(h)}_{i = 0} \psi_{i} \varepsilon_{T + h - i} - \sum^{K}_{h = 1} w_h\sum^{P(h)}_{i = h} \sum^{P(h)-i}_{j = 0} \sum^{P(h)-i-j}_{k = 0}  \widehat{\psi}_{i} \widehat{\pi}_j \psi _k \varepsilon_{T + h - i-j-k}\right) ^2 \right] \nonumber \\  
 & = E  \Biggr[\left(\sum^{K}_{h = 1} w_h\sum^{P(h)}_{i = 0} \psi_{i} \varepsilon_{T + h - i} \right) ^2 - 2 \sum^{K}_{h = 1} \sum^{K}_{h' = 1} w_h w_{h'} \sum^{P(h')}_{l = 0} \sum^{P(h)}_{i = h} \sum^{P(h)-i}_{j = 0} \sum^{P(h)-i-j}_{k = 0} \psi _l  \widehat{\psi}_{i} \widehat{\pi}_j \psi _k \varepsilon _{T+h'-l} \varepsilon_{T + h - i-j-k} \nonumber \\
 &+ \left( \sum^{K}_{h = 1}w_h \sum^{P}_{i = h} \sum^{P-i}_{j = 0} \sum^{P-i-j}_{k = 0}  \widehat{\psi}_{i} \widehat{\pi}_j \psi _k \varepsilon_{T + h - i-j-k}\right)^2  \Biggr]  = E  \Biggr[\left(\sum^{K}_{h = 1}\sum^{K}_{h' = 1} w_h w_{h'}\sum^{P(h)}_{i = 0}\sum^{P(h')}_{i' = 0} \psi_{i} \psi _{i'}\varepsilon_{T + h - i}\varepsilon_{T + h' - i'} \right) ^2 \Biggr]  \\
 &-2 \sigma ^2 \sum^{K}_{h = 1}\sum^{K}_{h' = 1} w_h w_{h'}\sum^{P(h)}_{l = 0} \sum^{P(h')}_{i = h'} \sum^{P(h)-i}_{j = 0} \sum^{P(h)-i-j}_{k = 0}  \psi_{l} \psi_k E \left[ \widehat{ \psi _i } \widehat{\pi }_j  \right] \delta _{h-l, h'-i-j-k}\nonumber \\
 &+E \Biggr[ \sum^{K}_{h = 1} \sum^{K}_{h' = 1} w_h w_{h'}\sum^{P(h)}_{i = h} \sum^{P(h)-i}_{j = 0} \sum^{P(h)-i-j}_{k = 0} \sum^{P(h')}_{i' = h'} \sum^{P(h')-i'}_{j' = 0} \sum^{P(h')-i'-j'}_{k' = 0} \psi _k \psi _{k'}  \widehat{ \psi} _i \widehat{\pi }_j  \widehat{\psi} _{i'} \widehat{\pi }_{j'} \varepsilon _{T+h-i-j-k} \varepsilon _{T+h'-i'-j'-k'}  \Biggr] \nonumber \\
 & = \sigma ^2 < \mathbf{w},(A+B+C) \mathbf{w}>,
 \end{align}
 where $A$, $B$, $C$ are the matrices with components given by
 \begin{align} 
\label{theor2:proof:23}
A_{hh'} = \sum^{P(h)}_{i = 0} \sum^{P(h')}_{i' = 0} \psi _i  \psi _{i'} \delta _{h-i, h'-i'},
\end{align}
\begin{align} 
\label{eq:theor423}
B_{hh'} = -2 \sum^{P(h)}_{l = 0} \sum^{P(h')}_{i = h'} \sum^{P(h)-i}_{j = 0} \sum^{P(h)-i-j}_{k = 0} \psi _l \psi _{k} E \left[ \widehat{ \psi }_i  \widehat{\pi} _j \right] \delta _{h-l, h'-i-j-k},
\end{align}
\begin{align} 
\label{eq:theor424}
C_{hh'} = \sum^{P(h)}_{i = h} \sum^{P(h)-i}_{j = 0} \sum^{P(h)-i-j}_{k = 0} \sum^{P(h')}_{i' = h'} \sum^{P(h')-i'}_{j' = 0} \sum^{P(h')-i'-j'}_{k' = 0} \psi _k \psi _{k'} E \left[ \widehat{ \psi }_i  \widehat{\pi} _j \widehat{\psi}_{i'} \widehat{\pi} _{j'} \right]  \delta _{h-i-j-k,h'-i'-j'-k'},
\end{align}
and $P(h) = T+h-1+r$, $P(h')=T+h'-1+r$.

\item[\bf{(iii)}]
We recall that by~(\ref{eq:theor41}), the forecast $ \widehat{X_{T + K}^\mathbf{w}}$ is given by
$$\widehat{X_{T + K}^\mathbf{w}} = \sum^{K}_{h' = 1} w_h\sum^{P(h)}_{i= h}\widehat{\psi} _i \tilde{\varepsilon} _{T+h-i}.$$
According to \eqref{eq:theor33}, both $ \widehat{\psi} _i$ and $ \widehat{\pi}_j $ can be asymptotically written as  
$$\widehat{\psi }_i = \psi_i + \dfrac{r_i}{\sqrt{T}}\quad \mbox{and} \quad \widehat{\pi }_j = \pi_j + \dfrac{t_j}{\sqrt{T}},$$
with $ r_i$ and $t_j$ Gaussian random variables of mean 0 and variances $(\Sigma _{\mathbf{\Xi}_{P}})_{i,i}$ and $(\Sigma _{\mathbf{\Xi}_{P}})_{j+P(K),j+P(K)}$, respectively. Consequently,  
\begin{align} 
& \widehat{ X_{T + K}^\mathbf{w}} = \sum^{K}_{h= 1} w_h\sum^{P(h)}_{i = h} \sum^{P(h)-i}_{j = 0} \sum^{P(h)-i-j}_{k = 0} \left(\psi_i + \dfrac{r_i}{\sqrt{T}} \right) \left( \pi_j + \dfrac{t_j}{\sqrt{T}}\right) \psi _k \varepsilon _{T+h-i-j-k}  \nonumber \\
& = \sum^{K}_{h = 1} w_h \sum^{P(h)}_{i = h} \sum^{P(h)-i}_{j = 0} \sum^{P(h)-i-j}_{k = 0} \left(\psi_i  \pi_j \psi _k+ \dfrac{\psi_i t_j  \psi _k}{\sqrt{T}} + \dfrac{r_i  \pi_j \psi _k}{\sqrt{T}} + \dfrac{r_i t_j  \psi _k}{T}\right) \varepsilon _{T+h-i-j-k}.\label{expression to simplify}
\end{align} 
We now recall that 
\begin{equation*}
\sum^{P(h)-i}_{j = 0} \sum^{P(h)-i-j}_{k = 0} \pi_j \psi _k \varepsilon _{T+h-i-j-k} = \varepsilon _{T+h-i},
\end{equation*}
and we eliminate in~(\ref{expression to simplify}) the term that decays as 1/T. We hence approximate $ \widehat{X_{T+K}^\mathbf{w}}$ as
\begin{align}
& \widehat{X_{T + K}^\mathbf{w}} \simeq\sum^{K}_{h= 1} w_h\sum^{P(h)}_{i = h} \left[ \left( \psi_i + \dfrac{r_i}{\sqrt{T}}  \right) \varepsilon _{T+h-i} + \sum^{P(h)-i}_{j = 0} \sum^{P(h)-i-j}_{k = 0} \dfrac{\psi _i \psi _k t _j }{\sqrt{T}} \varepsilon _{T+h-i-j-k}\right]   \nonumber \\
& = \sum^{K}_{h= 1} w_h\left[ \sum^{P(h)}_{i = h}  \widehat{\psi}_i \varepsilon _{T+h-i} + \sum^{P(h)}_{i = h}  \sum^{P(h)-i}_{j = 0} \sum^{P(h)-i-j}_{k = 0} \dfrac{\psi _i \psi _k t _j }{\sqrt{T}} \varepsilon _{T+h-i-j-k} \right] .  
\end{align}
Using this approximation we now compute the MSFE:
\begin{align}
 & \text{MSFE} \left( \widehat{X_{T + K}^\mathbf{w}} \right) = E \left[ \left( X _{T+K}^\mathbf{w} - \widehat{ X _{T+K}^\mathbf{w}} \right) ^2 \right]   \nonumber \\
 & =  E \left[ \left(\sum^{K}_{h= 1} w_h \left[ \sum^{P(h)}_{i = 0} \psi_{i} \varepsilon_{T + h - i} - \sum^{P(h)}_{i = h}  \widehat{\psi}_i \varepsilon _{T+h-i} - \sum^{P(h)}_{i = h}  \sum^{P(h)-i}_{j = 0} \sum^{P(h)-i-j}_{k = 0} \dfrac{\psi _i \psi _k t _j }{\sqrt{T}} \varepsilon _{T+h-i-j-k}  \right] \right) ^2 \right] \nonumber \\
 & = E \left[ \left(\sum^{K}_{h= 1} w_h \left[ \sum^{h-1}_{i = 0} \psi_{i} \varepsilon_{T + h - i} + \sum^{P(h)}_{i = h}  \left( \psi _i - \widehat{\psi}_i \right) \varepsilon _{T+h-i} - \sum^{P(h)}_{i = h}  \sum^{P(h)-i}_{j = 0} \sum^{P(h)-i-j}_{k = 0} \dfrac{\psi _i \psi _k t _j }{\sqrt{T}} \varepsilon _{T+h-i-j-k} \right]  \right) ^2 \right] \nonumber \\
 & = E \left[ \sum^{K}_{h= 1} \sum^{K}_{h'= 1} w_h w_{h'}  \sum^{h-1}_{i = 0} \sum^{h'-1}_{i' = 0} \psi_{i} \psi _{i'} \varepsilon_{T + h - i}\varepsilon_{T + h' - i'} \right]  \nonumber  \\
 & + E \left[ \sum^{K}_{h= 1} \sum^{K}_{h'= 1} w_h w_{h'} \sum^{P(h)}_{i = h} \sum^{P(h')}_{i' = h'}  \left( \psi _i - \widehat{\psi}_i \right) \left( \psi _{i'} - \widehat{\psi}_{i'} \right) \varepsilon _{T+h-i} \varepsilon _{T+h'-i'} \right]  \nonumber \\
 & + E \left[  \sum^{K}_{h= 1} \sum^{K}_{h'= 1} w_h w_{h'} \sum^{P(h)}_{i = h}  \sum^{P(h)-i}_{j = 0} \sum^{P(h)-i-j}_{k = 0} \sum^{P(h')}_{i' = h'}  \sum^{P(h')-i'}_{j' = 0} \sum^{P(h')-i'-j'}_{k' = 0} \dfrac{\psi _i \psi _k \psi _{i'} \psi _{k'} t _j t_{j'}}{{T}} \varepsilon _{T+h-i-j-k} \varepsilon _{T+h'-i'-j'-k'} \right] \nonumber \\
 & + 2 E \left[  \sum^{K}_{h= 1} \sum^{K}_{h'= 1} w_h w_{h'} \sum^{P(h)}_{i = h} \sum^{P(h')}_{i' = h'}  \sum^{P(h')-i'}_{j' = 0} \sum^{P(h')-i'-j'}_{k' = 0} \psi _{i'} \psi _{k'} \psi _{i'} \left( \widehat{\psi }_i - \psi _i \right) \left( \widehat{\pi}_{j'}  - \pi _{j'} \right) \varepsilon _{T+h-i} \varepsilon _{T+h'-i'-j'-k'} \right] \nonumber.
\end{align}
Using Lemma~\ref{eq:theor33}, this expression can be asymptotically approximated by:  
\begin{equation*} 
\text{MSFE} \left( \widehat{ X_{T + K}^\mathbf{w} }\right) = \sigma ^2 < \mathbf{w}, \left( A^{char} + D + F + G \right) \mathbf{w}>,
\end{equation*}
where $A^{char}$, $D$, $F$, $G$ are matrices whose components are given by
$$A_{hh'}^{char}:= \sum^{h-1}_{i = 0}  \sum^{h'-1}_{i' = 0} \psi _i \psi _{i'} \delta _{h-i,h'-i'},$$
$$D_{hh'}=\dfrac{1}{T} \sum^{P(h)}_{i = h} \sum^{P(h')}_{i' = h'} \left( \Sigma _ {\mathbf{\Xi}_{P}} \right)_{i,i'}   \delta _{h-i,h'-i'},$$
$$F_{hh'} = \dfrac{2}{T} \sum^{P(h)}_{i=h} \sum^{P(h')}_{i' = h'} \sum^{P(h')-i'}_{j' = 0} \sum^{P(h')-i'-j'}_{k' = 0} \psi _{i'} \psi _{k'} \left( \Sigma _{{\Xi}_P}\right)_{i, P(K)+j'}  \delta _{h-i, h'-i'-j'-k'},$$
$$\!\!\!\!\!\!\!\!\!\!\!\!\!\!\!\!G_{hh'} = \dfrac{1}{T}\sum^{P(h)}_{i = h} \sum^{P(h)-i}_{j = 0} \sum^{P(h)-i-j}_{k = 0} \sum^{P(h')}_{i' = h'} \sum^{P(h')-i'}_{j' = 0} \sum^{P(h')-i'-j'}_{k' = 0} \psi _i \psi _{k} \psi_{i'} \psi _{k'}   \left( \Sigma _ {\mathbf{\Xi}_{P} }\right)_{j+P(K),j'+P(K)}    \delta _{h-i-j-k,h'-i'-j'-k'}.$$
\end{enumerate}

\addcontentsline{toc}{section}{Bibliography}
\bibliographystyle{alpha}
\bibliography{/Users/JP17/JPO_synch/BiblioData/Library-GOLibrary}
\end{document}